\documentclass{article}

\usepackage{amsmath,amsfonts,epsfig,enumerate}
\usepackage{mathptmx}
\usepackage{amsmath,amsfonts,epsfig,enumerate}
\usepackage{mathptmx}

\usepackage{amsmath,amsfonts,latexsym,amssymb}
\usepackage{mathrsfs}
\usepackage{verbatim}
\usepackage[utf8]{inputenc}
\usepackage[T1]{fontenc}
\usepackage{ae,aecompl}
\usepackage{braket}
\usepackage{color}
\usepackage{euscript}

\newcommand{\widemargins}{
\setlength{\textwidth}{5.8in}
\setlength{\oddsidemargin}{0.25in}
\setlength{\evensidemargin}{0.25in}
}

\renewcommand{\bold}[1]{\medskip \noindent {\bf \boldmath #1
                        }\nopagebreak[4]}

\newcommand{\grad}{\nabla}
\newcommand{\csch}{\operatorname{csch}}

\newcommand{\sech}{\operatorname{sech}}

\newcommand{\Teich}{\operatorname{Teich}}

\newtheorem{theorem}{Theorem}[section]
\newtheorem{prop}[theorem]{Proposition}
\newtheorem{lemma}[theorem]{Lemma}
\newtheorem{cor}[theorem]{Corollary}

\newtheorem{corollary}[theorem]{Corollary}

\renewcommand{\hbar}{\bar{{\mathbb H}}^3}

\newcommand{\R}{\mathbb R}

\newcommand{\N}{\mathbb N}

\newcommand{\Hp}{{\mathbb H}^2}

\newcommand{\T}{{\operatorname{Teich}}}

\def\lb{6.57252}
\def\ub{6.65603}
\def\lbs{9.29495}
\def\ubs{9.41305}

\def\sys{{\operatorname{sys}}}
\def\eproof{$\Box$ \medskip}

\newcommand{\ssm}{\smallsetminus}

\widemargins

\begin{document}
\title{\bf \Large Strata separation for the Weil-Petersson completion and gradient estimates for length functions} \author{Martin
   Bridgeman\thanks{M. Bridgeman's research was supported by NSF grant DMS-1500545.} \  and
  Kenneth Bromberg\thanks{K. Bromberg's research supported by NSF grant DMS-1906095.}}

\date{\today}

\maketitle

\begin{abstract}
\noindent
In general,  it is difficult to measure distances in the Weil-Petersson metric on Teichm\"uller space.
Here we consider the distance between strata in the Weil-Petersson completion of Teichm\"uller space of a surface of finite type. Wolpert showed that for strata whose closures do not intersect, there is a definite separation independent of the topology of the surface. 
We prove that the optimal value for this minimal separation is a constant $\delta_{1,1}$ and show that  it is realized exactly by strata whose nodes intersect once. We also give a nearly sharp estimate for $\delta_{1,1}$
 and give a lower bound on the size of the gap between $\delta_{1,1}$ and the other distances. A major component of the paper is an effective version of Wolpert's upper bound on $ \langle \grad \ell_\alpha,\grad \ell_\beta \rangle$, the inner product of the Weil-Petersson gradient of length functions. We further bound the distance to the boundary of Teichm\"uller space of a hyperbolic surface in terms of the length of the systole of the surface. We also obtain new lower bounds on the systole for the Weil-Petersson metric on the moduli space of a punctured torus. \end{abstract}

\section{Strata separation}
There are several natural quantities associated to the Weil-Petersson metric on Teichm\"uller and moduli space. One is the length of closed geodesics on moduli space or, equivalently, the translation length of pseudo-Anosovs on Teichm\"uller space. Another is the distance between {\em strata} on the boundary of Teichm\"uller space. Boundary strata are determined by a multi-curve on the underlying  surface and two strata will have intersecting closures if and only if the associated multi-curves have positive intersection. Wolpert has shown that there is a definite separation (independent of the surface) between two strata whose closures do not intersect. The key tool in the proof of this theorem are upper bounds on the gradients of length functions. In this paper we will improve on Wolpert's gradient estimates and use this to show that, as expected, the minimal distance is realized when the multi-curves intersect exactly once. We will also see that nearly sharp bounds on  this distance follow easily for our gradient estimates.

We begin with some setup before stating our results more precisely.
Let $S$ be hyperbolic surface of finite type and $\Teich(S)$ the associated Teichm\"uller space. We let   $\overline{\Teich(S)}$ be the completion with respect to the  Weil-Petersson metric. 

There is a natural stratification of $\overline{\Teich(S)}$ which can be described via length functions.

Given a closed curve (or multi-curve) $\alpha$ in $S$ we have the length function $\ell_\alpha\colon \Teich(S)\rightarrow (0,\infty)$ given by letting $\ell_\alpha(X)$ be  the length of the geodesic representative of $\alpha$ in $X$. Then $\ell_\alpha$ extends to a continuous function $\ell_\alpha\colon \overline{\Teich(S)}\rightarrow [0,\infty]$. Given a multi-curve $\tau$ on $S$, we define the associated stratum 
$$\mathcal S_\tau(S) = \{ X \in \overline{\Teich(S)} \mbox{ such that } \ell_{\alpha}(X) = 0 \mbox{ if and only if } \alpha \subseteq \tau\}.$$ 
Points in $\mathcal S_\tau(S)$ are {\em noded} hyperbolic structures on $S$ where the multi-curve $\alpha$ is the set of nodes.

We note that if $\sigma \subseteq \tau$ then $\mathcal S_\tau(S) \subseteq \overline{\mathcal S_\sigma(S)}$ and  it follows easily that   $i(\sigma,\tau) = 0$ if and only if $d_{\rm WP}(\mathcal S_\sigma(S),\mathcal S_\tau(S)) = 0$. 

Wolpert proved the following:
\begin{theorem}[{Wolpert Strata Separation, \cite{Wolpert:compl}}]
There is a universal constant $\delta_0 > 0$ such that if $\mathcal S_\sigma(S),\mathcal S_\tau(S)$ are two strata with geometric intersection number $i(\sigma,\tau) \neq 0$ then $d_{\rm WP}(\mathcal S_\sigma(S),\mathcal S_\tau(S)) \geq \delta_0$.
\label{wolpert:strata}
\end{theorem}

Wolpert does not give an explicit value for the constant $\delta_0$. We will give the optimal value for $\delta_0$.

We let $T$ be a punctured torus and $\alpha,\beta$ two curves on $T$ with $i(\alpha,\beta)=1$. Observe that there is an element of the mapping class group (i.e. an isometry of $\overline{\Teich(T)}$) that takes any other pair of curves on $T$ that intersect once to $\alpha$ and $\beta$ so the constant
$$\delta_{1,1} = d_{\rm WP}(\mathcal S_\alpha(T),\mathcal S_\beta(T))$$
is well defined. An elementary application of Riera's formula (see Lemma \ref{dbound})  shows that $\lb < \delta_{1,1} < \ub.$

Using estimates on the Weil-Petersson gradient of  length functions  along with Wolpert's description of the 
Alexandrov tangent cone for the Weil-Petersson completion, we prove that the optimal value for Wolpert's constant $\delta_0$ is exactly $\delta_{1,1}$. More precisely:
  \begin{theorem}
  Let  $\mathcal S_\sigma(S),\mathcal S_\tau(S)$ be two strata in $\Teich(S)$. Then one of the following holds;
  \begin{enumerate}
\item  $i(\sigma,\tau) = 0$ and  $d_{\rm WP}(\mathcal S_\sigma(S),\mathcal S_\tau(S))  = 0.$
\item $i(\sigma,\tau) = 1$ and  $d_{\rm WP}(\mathcal S_\sigma(S),\mathcal S_\tau(S))  = \delta_{1,1}.$
\item $i(\sigma,\tau) > 1$ and  $d_{\rm WP}(\mathcal S_\sigma(S),\mathcal S_\tau(S)) \geq 7.61138.$ 
 \end{enumerate}
 \label{smallsep}
\end{theorem}
We note that it is not hard to see that the set of distances between strata (even for the punctured torus) is not a discrete set and Wolpert's original theorem does not give that the constant $\delta_0$ is attained.

If $S$ is a punctured sphere then intersecting curves intersect at least twice and this setting needs a slightly separate analysis. See section \ref{nps}.

Another application is relating the distance of a point  in $\mathcal \T(S)$ from the boundary $\partial\overline{\T(S)}$ to the length of its systole. Given $X \in \mathcal \T(S)$ we let $\ell_{sys}(X)$ be the length of the systole of $X$, i.e. the minimum length of a geodesic on $X$. We prove

\begin{theorem}
There exists an explicit continuous function $c:(0,\infty) \rightarrow (0,1)$ such that if $S$ is a surface of finite type and $X \in \T(S)$ then
$$\sqrt{\frac{2}{\pi}}\leq c(\ell_{sys}(X)) \leq \frac{d_{\rm WP}(X,\partial\overline{\T(S)})}{\sqrt{2\pi\ell_{sys}(X)}} \leq 1.$$
Furthermore $\lim_{t\rightarrow 0}c(t) = \lim_{t\rightarrow \infty}c(t)=1$.
\end{theorem}

The in-radius of $\T(S)$ is the radius of the largest embedded metric ball in  $\T(S)$ (see \cite{Brock:Bromberg:vol} and \cite{Wu:inradius}). Specifically
$$\mbox{\rm InRad}(\T(S)) = \max_X d_{\rm WP}(X,\partial\overline{\T(S)}).$$
If we let ${\rm sys}(S) = \max_{X \in \T(S)} \ell_{sys}(X)$ then the above theorem gives the following immediate corollary.

\begin{corollary} With $c$ the same as above
$$c({\rm sys}(S)) \leq \frac{\mbox{\rm InRad}(\T(S))}{\sqrt{2\pi{\rm sys}(S)}} \leq 1.$$
Furthermore if $S_{g,n}$ is the surface of type $g,n$ then
$$\lim_{g\rightarrow \infty} \frac{\mbox{\rm InRad}(\T(S_{g,n}))}{\sqrt{2\pi{\rm sys}(S_{g,n})}} = 1.$$
\end{corollary}

\subsubsection*{Gradient estimates}
Riera gave a beautiful formula for the inner product of the Weil-Petersson gradient of length functions $\ell_\alpha$ and $\ell_\beta$
(see Theorem \ref{Riera}). Using this formula Wolpert obtained the following estimate:
\begin{theorem}[{Wolpert, \cite{Wolpert:behavior}}]
Let $\ell_\alpha, \ell_\beta$ be geodesic length functions for simple disjoint  curves $\alpha,\beta$.  Then
$$\frac{2}{\pi}\ell_\alpha(X)\delta^\alpha_\beta  \leq \langle \grad \ell_\alpha,\grad \ell_\beta \rangle \leq \frac{2}{\pi}\ell_\alpha(X)\delta^\alpha_\beta +O(\ell_\alpha(X)^2\ell_\beta(X)^2)$$
where $\delta^\alpha_\beta$ is the Kronecker delta function and where for $\ell > 0$   the term $O(\ell_\alpha(X)^2l_\beta(X)^2)$ is uniform  for $\ell_\alpha(X),\ell_\beta(X) < \ell$.
\label{wgrad}
\end{theorem}
The lower bound follows directly from Riera's formula. Following the same basic strategy of Wolpert's proof we obtain an upper bound on the inner product by an explicit elementary function. As in Wolpert's bound this function will decay quadratically in both $\ell_\alpha(X)$ and $\ell_\beta(X)$ as the lengths approach zero but for large lengths it grows exponentially:
\begin{theorem}
Let $\ell_\alpha, \ell_\beta$ be geodesic length functions for simple disjoint  curves $\alpha, \beta$ with $ \ell_\alpha(X)\leq \ell_\beta(X)$. Then
$$\frac{2}{\pi}\ell_\alpha(X)\delta^\alpha_\beta  \leq \langle \grad \ell_\alpha,\grad \ell_\beta \rangle \leq \frac{2}{\pi}\ell_\alpha(X)\delta^\alpha_\beta +\frac{8}{3\pi^2}\ell_\alpha(X)\sinh\left(\ell_\alpha(X)/2\right)\sinh^2\left(\ell_\beta(X)/2\right)$$
where $\delta^\alpha_\beta$ is the Kronecker delta function.
\label{grad}
\end{theorem}
We note that the bound here is asymptotically optimal for small lengths but not when the length is large. One can obtain a better bound by an elementary (but complicated) function that has better asymptotics for large lengths (see Proposition \ref{elementary_bound}). At the end of Section 2 there is a further discussion on the accuracy of our bounds.

\subsection*{Notation}
In using decimals approximations the expression $a \simeq a_0.a_1a_2a_3\ldots a_n$  where $a_0 \in \N_0$ and $a_i \in \{0,1,\ldots,9\}$  means that this  is the  first $n$ decimal places of $a$.

\subsubsection*{Acknowledgements}
The authors would like to thank Jeffrey Brock,  Scott Wolpert and Yunhui Wu  for helpful conversations on this project.

\section{Bounding the gradient}

\subsubsection*{Riera's formula}
The main ingredient in Wolpert's bound is the following formula of Riera for the Weil-Petersson inner product of length functions $\ell_\alpha$ and $\ell_\beta$.

\begin{theorem}[{Riera,  \cite{Riera:formula}}]
For $X = \Hp/\Gamma \in \Teich(S)$, let
$A,B \in \Gamma$  correspond to $\alpha,\beta$ with $A,B$ having axes $a,b$. For $C \in \langle  A\rangle \setminus \Gamma/\langle B\rangle$, if $a,C(b)$ intersect  let $u(C) = \cos(a,C(b))$  the cosine of the angle of intersection and otherwise let $u(C) = cosh(d(a,C(b))$. 
Then
$$\langle \grad \ell_\alpha,\grad \ell_\beta \rangle_X = \frac{2}{\pi}\left(\ell_\alpha(X)\delta^\alpha_\beta +\sum_{C \in \langle A\rangle\setminus \Gamma / \langle B\rangle} R(u(C))\right)$$
where
$$R(u) = u\log\left|\frac{u+1}{u-1}\right| -2 .$$
\label{Riera}\end{theorem}

Before starting on the main estimate of the paper we use Riera's formula to bound the distance between strata in a simple, but important case. We note that if $u>1$, $R(u)>0$ so if the curves $\alpha$ and $\beta$ are disjoint (or equal) then the inner product of $\nabla \ell_\alpha$ and $\nabla \ell_\beta$ will be positive.

\begin{prop}
The constant $\delta_{1,1}$ has the following bounds:
$$\sqrt{\pi}\int_0^{2\sinh^{-1}(1)} \frac{dt}{\sqrt{\sinh(t/2)}} \le \delta_{1,1} \le 4\sqrt{\pi\sinh^{-1}(1)}$$
In particular, numerical estimates give $\delta_{1,1} \in (\lb,\ub).$
\label{dbound}\end{prop}

{\bf Proof:} 
Let $\alpha$ and $\beta$ be curves on the punctured torus $T$ that intersect once. There is orientation reversing involution $\iota \colon T\to T$ that fixes both $\alpha$ and $\beta$ (as homotopy classes). This involution induces an isometric involution $\iota_*\colon \Teich(T) \to \Teich(T)$ under which both $\ell_\alpha$ and $\ell_\beta$ (and therefore their gradients) are invariant. (In fact there are two such involutions $\iota$ but they both induce the same map on $\Teich(T)$.) We construct a path $X_t$ in $\Teich(T)$ from $\mathcal S_\alpha$ to $\mathcal S_\beta$ (which are both single points) that is the fixed point set of $\iota_*$. This implies that $X_t$ is the unique geodesic from $\mathcal S_\alpha$ to $\mathcal S_\beta$ and that both $\nabla \ell_\alpha$ and $\nabla \ell_\beta$ are tangent to it.

Here is a description of the path: Let $R_t$ be a family of ideal quadrilaterals where the two shortest geodesics connecting opposite sides intersect orthogonally and one of these sides has length $t$ and the other $s$. A direct calculation shows that
$$\sinh(t/2)\sinh(s/2) = 1.$$
The tori $X_t$ are defined by identifying  the opposite sides of $R_t$. The involution $\iota$ is induced by reflecting $R_t$ along the horizontal geodesic which induces an isometry of $X_t$ to itself. For any other torus $Y\in \Teich(T)$ the angle between $\alpha$ and $\beta$ will be some $\theta \neq \pi/2$ while the angle $\iota(\alpha)$ and $\iota(\beta)$ will be $\pi - \theta$ so $\iota(Y) \neq Y$. Therefore $X_t$ is the fixed point set of $\iota_*$.

Note that $\ell_\alpha(X_t) = t$ and $\ell_\beta(X_t) = s$ so on $X_t$ the relationship between $t$ and $s$ gives
$$\sinh(\ell_\alpha(X_t)/2)\sinh(\ell_\beta(X_t)/2) =1.$$
This and the Riera Formula will allow us to get good bounds on the gradients. In particular, given that the gradients are tangent to $X_t$, after differentiating we have
$$\nabla \ell_\alpha(X_t) = -\sinh(\ell_\alpha(X_t)/2)\nabla \ell_\beta(X_t).$$
Applying Riera's formula to the inner product of $\nabla \ell_\alpha$ with itself we have
$$\|\nabla \ell_\alpha(X_t)\|^2 \ge \frac2\pi \ell_\alpha(X_t)=\frac{2t}\pi$$
as all of the terms in the sum are positive. If we take the the inner produce of $\nabla\ell_\alpha$ and $\nabla\ell_\beta$ we have
$$\langle \nabla\ell_\alpha,\nabla \ell_\beta\rangle \ge -\frac4\pi$$
as the only non-positive term comes from lift of  $\beta$ that intersects the lift of $\alpha$ in the double coset. As the two gradients are tangent but in opposite directions we also have
$$\langle \nabla\ell_\alpha,\nabla\ell_\beta\rangle = - \|\nabla \ell_\alpha\| \cdot \|\nabla\ell_\beta\|.$$
Combining with our previous relationship on the gradients we have
$$\|\nabla\ell_\alpha(X_t)\|^2 \le \frac4\pi \sinh(\ell_\alpha(X_t)/2) = \frac4\pi \sinh(t/2).$$

Choose $t_0 = 2\sinh^{-1}(1)$. Then $\ell_\alpha(X_{t_0}) = \ell_\beta(X_{t_0})$ and by symmetry the length of the paths $X_{(0, t_0]}$ and $X_{[t_0, \infty)}$ are equal. We will use the above bounds on gradients to bound the length of the former.

As the tangent vector $\dot X_t$ is parallel to $\nabla \ell_\alpha(X_t)$ after differentiating the formula $\ell_\alpha(X_t) = t$ we have
$$\|\dot X_t \| \cdot \|\nabla \ell_\alpha(X_t)\|  = 1$$
and therefore
$$\operatorname{Length}\left(X_{(0,t_0]}\right) = \int_0^{t_0} \frac1{\|\nabla \ell_\alpha(X_t)\|} dt.$$
Applying our estimates on $\|\nabla \ell_\alpha(X_t)\|$ we have
$$\int_0^{t_0} \frac1{\sqrt{\frac4\pi \sinh(t/2)}} dt \le \operatorname{Length}\left(X_{(0,t_0]}\right) \le\int_0^{t_0} \frac1{\sqrt{\frac{2t}\pi}} dt.$$
As $\delta_{1,1} = 2\operatorname{Length}\left(X_{(0,t_0]}\right)$, the result follows. \eproof

{\em Remark:} The first bounds on $\delta_{1,1}$ we given in \cite{Brock:Bromberg:vol} where it was shown that
$$0.9744... < \delta_{1,1}<25.8496...$$
using bounds on volumes hyperbolic 3-manifolds. The method here allows one to estimate $\delta_{1,1}$ to any degree of accuracy.  As $\pi_1(T) =\langle \alpha,\beta\rangle$  we can enumerate the double cosets in Riera's formula for both $\|\grad \ell_\alpha\|^2$ and  $\langle \nabla\ell_\alpha,\nabla \ell_\beta\rangle$ in terms of words in $A$ and $B$. These enumerations give  distance functions $u_i(t)$   and $v_i(t)$  so that for any $m,n \in \N_0$
$$\frac{2}{\pi}\sinh(t/2)\left( 2 -  \sum_{i=1}^m R\left(v_i(t)\right) \right)\leq \|\grad \ell_\alpha(X_t)\|^2 \leq \frac{2}{\pi}\left(t + \sum_{i=1}^n R(u_i(t))\right).$$
 In particular taking the double cosets $C_n = \langle A\rangle\setminus B^n / \langle B\rangle$ then  $u(C_n) = \cosh(n\ell_\beta)$.  Similarly, we observe that the 4 double cosets $C_{\pm,\pm} = \langle A\rangle\setminus B^{\pm}A^{\pm} / \langle B\rangle$ give  $u(C_{\pm,\pm}) =\sinh(\ell_\alpha)\sinh(\ell_\beta).$ 
Using these upper and lower bounds, we numerically integrate to obtain
$$6.59576 \leq \delta_{1,1} \leq 6.63283.$$

In this example the upper bounds on $\|\nabla \ell_\alpha\|$ are obtained by exploiting the extra symmetry in this setting. To bound $\|\nabla \ell_\alpha\|$ in a more general setting (which we use to bound distances in $\Teich(S)$) we need to bound the sum in Riera's formula directly.

\subsubsection*{Strategy}
We briefly describe the strategy for the proof of Theorem \ref{grad}. The function $R(\cosh(t))$ can be approximated by $ae^{-2t}$. To bound the sum in Riera's formula we compare it to the integral of the function $e^{-2d(\alpha, z)}$ on the annular cover $A_\alpha$ of $X$ associated to $\alpha$ where $d(\alpha, z)$ is the distance be a point $z\in A_\alpha$ and the core geodesic. The integral over the annulus is a straightforward calculation. To compare it to the sum we decompose the annulus into the $r$-neighborhoods $N(h_i, r)$ of the lifts $h_i$ of $\beta$ to $A_\alpha$ where $r$ is an explicit constant given by the collar lemma and then compare the average value of $e^{-2d(\alpha, z)}$ on $N(h_i, r)$ to $e^{-2d_i(X)}$.

While the overall strategy of the proof is the same as Wolpert's, our estimates within the proof are different. For example Wolpert only estimates the average of $e^{-2d(\alpha,z)}$ on disks rather than over the neighborhoods $N(h_i, r)$.

\subsubsection*{Preliminary estimates}
Before proving the theorem we need to approximate $R$ and implement our averaging estimate. We begin with the former.
\begin{lemma}
The function 
$$a(t) = e^{2t}R(\cosh t)$$
is monotonically decreasing with 
$$\underset{t\to\infty}{\lim}\  a(t) = \frac{8}{3}.$$
Furthermore
$$a(t) \leq \frac{8}{3}-2\log(1-e^{-2t}).$$
\label{eRiera}
\end{lemma}

{\bf Proof:}
We have by \cite{Riera:formula} that for $s > 1$ 
$$R(s) = s\log\left(\frac{s+1}{s-1}\right) -2 =  \frac{2}{3s^2} + \frac{2}{5s^4} +\frac{2}{7s^6}+\ldots.$$
Note that if we replace $R$ by its series above, the individual terms of $e^{2t}R(\cosh(t))$ are not each monotonically decreasing. To prove the lemma we need a different expansion of $a(t)$. Let $u = e^{-t}$ and consider
$$\hat a(u) =  u^{-2}R\left(\frac{u+1/u}{2}\right).$$
We have
\begin{eqnarray*}
R\left(\frac{u+1/u}{2}\right)  &=&  \left(\frac{u+1/u}{2}\right)\log\left(\frac{\left(\frac{u+1/u}{2}\right)+1}{\left(\frac{u+1/u}{2}\right)-1}\right) -2\\
&=&  \left(\frac{u+1/u}{2}\right)\log\left(\frac{u^2+2u+1}{u^2-2u+1}\right) -2\\
&=& (u+1/u) \log\left(\frac{1+u}{1-u}\right) -2 \\
 &=& (u+1/u)\left(2u+ \frac{2u^3}{3} + \frac{2u^5}{5}+\ldots\right) -2\\
&=& \sum_{n=1}^\infty \left(\frac{2}{2n-1}+\frac{2}{2n+1}\right)u^{2n} = \sum_{n=1}^\infty \frac{8n}{(2n-1)(2n+1)}u^{2n}\\
\end{eqnarray*}
Therefore 
$$\hat a(u) = \sum_{n=0}^\infty \frac{8(n+1)}{(2n+1)(2n+3)}u^{2n}.$$
From the expansion, it follows that $\hat a(u)$ is monotonically increasing on $[0,1)$ and therefore $a(t) = \hat a\left(e^{-t}\right)$ is monotonically decreasing on $(0,\infty)$ and
$$ \underset{t\to\infty}{\lim}\  a(t) = \hat a(0) = \frac 8 3.$$ 
To obtain the upper bound, we have
\begin{eqnarray*}
\hat a(u) &=&   \frac{8}{3} + \sum_{n=1}^\infty \frac{8(n+1)u^{2n}}{(2n+1)(2n+3)} \\
&=& \frac{8}{3} + 2\sum_{n=1}^\infty \left(\frac{2}{2n+1}\right)\left(\frac{2n+2}{2n+3}\right)u^{2n}\\
 &\leq&  
\frac{8}{3} + 2\sum_{n=1}^\infty \frac{u^{2n}}{n} \\
&=& \frac{8}{3} - 2\log(1-u^2).
\end{eqnarray*}
\eproof

Let $d$ denote distance in the hyperbolic plane $\Hp$ and $dA$ the hyperbolic area form.
We will use the following lemma to estimate the integral of $e^{-d(\alpha,z)}$ over $N(h_i, r)$.
\begin{lemma}
Let $g, h$ be disjoint geodesics with $d(g,h) > r$ and let $N(h,r)$ be the $r$ neighborhood of $h$. Then
$$e^{2d(g,h)}  \int_{N(h,r)} e^{-2d(g,w)} dA \geq 2\tan^{-1}(\sinh(r))\cosh^2(r) +2\sinh(r).$$
Furthermore if $d(g,h_n) \rightarrow \infty$ then
$$\lim_{n\rightarrow \infty}\left( e^{2d(g,h_n)}\int_{N(h_n,r)} e^{-2d(g,w)} dA\right) = 2\tan^{-1}(\sinh(r))\cosh^2(r) +2\sinh(r).$$
\label{mean}
\end{lemma}
\noindent{\bf Proof:}
We first make a general observation. We consider the triple $(E, p , g)$ where $E$ is a Borel set in 
$\Hp$, $p \in E$ and $g$ is a geodesic such that $E$ is entirely on one side of $g$. Note that if $\mathfrak h$ is a horocycle tangent to $g$ that is on the other side of $E$ then $d(q, \mathfrak h) \ge d(q, g)$ for all $q\in E$. Therefore
$$\int_E e^{-2d(g,w)}dA \ge \int_E e^{-2d(\mathfrak h, w)} dA.$$
We can estimate the integral on the right by working in the half space model for $\Hp$ and normalizing so that $p=i$, $g$ intersects the imaginary axis at $y_0>1$ and $\mathfrak h$ is the horizontal line at height $y_0$. Then for $w = (x,y) \in E$ we have
$$d(\mathfrak h, w) = \log\left(\frac{y_0}y\right)$$
and
$$\int_E e^{-2d(\mathfrak h, w)}dA =  \int_{E_0} \frac{y^2}{y_0^2}\cdot \frac{dxdy}{y^2}= \frac{A(E)}{y_0^2} = e^{-2d(g,p)}A(E).$$

Let $g_n$ be a sequence of geodesics such that in the normalized picture $g_n$ intersect at height $y_n$ with $\lim_{n\rightarrow \infty} y_n = \infty$. Let $\mathfrak h_n$ be the horocycle for $y = y_n$. Then for $f_n(w) = d(w, \mathfrak h_n)-d(w, g_n)$  we have $f_n \rightarrow 0$ uniformly on compact subsets of $\Hp$. Therefore
$$\lim_{n\rightarrow \infty} \left(e^{2d(g_n,p)}\int_E e^{-2d(g_n, w)}dA\right)  =  \lim_{n\rightarrow \infty} \left(e^{2d(g_n,p)}\int_E e^{-2d(\mathfrak h_n, w)}dA\right)  = A(E).$$

We now apply this to the geodesics $g,h$. We consider the triple $(N(h,r), p, g)$ where $p$ is the nearest point on $h$ to $g$. Then from above
$$\int_{N(h,r)} e^{-2d(g,w)} dA \ge e^{-2d(g,h)}A(r)$$
where $A(r)$ is the Euclidean area of $N(h,r)$ when $h$ is the semicircle of radius $1$ about $0$.

To calculate $A(r)$, we do some basic calculus. We define $\phi$ to be the angle between the boundary of $N(h,r)$  and the geodesic $h$. Reflecting the bottom boundary component of $N(h,r)$ in the $x$-axis, the upper boundary component and the reflected bottom make a Euclidean circle of radius $R$ with $R\cos(\phi) = 1$. 
We then consider a Euclidean circle C of radius $R$ about the origin and let $I(t)$ be the area  between the vertical line $x= t$ and $C$. Then
$$I(t) = 2\int_t^R \sqrt{R^2-x^2}dx.$$
We observe that $A(r) = \pi R^2-2I\left(\sqrt{R^2-1}\right)$. Substituting $x= R\sin\theta$ we have
$$I\left(\sqrt{R^2-1}\right) =  2R^2\int_{\phi}^{\pi/2}\cos^2\theta d\theta = R^2\left(\frac{\pi}{2}-\phi-\frac{1}{2}\sin(2\phi)\right).$$
Thus
$$A(r) = R^2(2\phi+\sin(2\phi)) = \frac{2\phi+2\sin(\phi)\cos(\phi)}{\cos^2(\phi)}.$$
By elementary hyperbolic geometry $\cosh(r) = \sec(\phi), \sinh(r) = \tan(\phi)$ and $ \tanh(r) = \sin(\phi)$. Therefore
$$A(r) = 2\tan^{-1}(\sinh(r))\cosh^2(r) +2\sinh(r).$$
\eproof

 Using the above lemmas we prove the following proposition.
 \begin{prop}\label{elementary_bound}
 Let $\ell_\alpha, \ell_\beta$ be geodesic length functions for $\alpha,\beta$ simple and disjoint. Then
 $$\frac{2}{\pi}\ell_\alpha(X)\delta^\alpha_\beta  \leq \langle \grad \ell_\alpha,\grad \ell_\beta \rangle \leq \frac{2}{\pi}\ell_\alpha(X)\left(\delta^\alpha_\beta +F(\ell_\alpha, \ell_\beta)\right)$$
where $F$ is an explicit elementary function.
\end{prop}
  
{\bf Proof:}
We let $A_\alpha$ be the annular cover corresponding to geodesic $\alpha$ in $X$. We let $g$ be the core geodesic and $h_i$ an enumeration of the  lifts of $\beta$ in $A_\alpha$. We further let $t_i$ be the distance from $g$ and $h_i$.  Then by \cite{Riera:formula} we have: 
$$\langle \grad \ell_\alpha,\grad \ell_\beta \rangle = \frac{2}{\pi}\left(\ell_\alpha\delta^\alpha_\beta + \sum_{i} R(\cosh(t_i))\right).$$

The lower bound on $\langle \grad \ell_\alpha,\grad \ell_\beta \rangle$ then follows as $R(t) > 0$ for $t> 1$. 
We let $T$ be the minimum distance between $\alpha$ and $\beta$ and $r,s > 0$ be such that the $r$-neighborhood $\alpha$ and the $s$ neighborhood of $\beta$ are both embedded and disjoint. In particular $T \geq r+s$. Also by the collar lemma, $\sinh(r) \geq 1/\sinh(\ell_\alpha/2)$ and $\sinh(s) \geq 1/\sinh(\ell_\beta/2)$.  

As $d_i(X) \geq T$ for all $i$, by the  Lemma \ref{eRiera} 
$$\sum_{i} R(\cosh(t_i)) \leq a(T)\sum_{i=1}^\infty   e^{-2t_i} \leq a(r+s)\sum_{i=1}^\infty   e^{-2t_i}$$
We now bound the expression on the right. 

Define  $N(h_i,s)$ to be the $s$-neighborhood of  $h_i$ and $N(g,r)$ to be the $r$-neighborhood of $g$. Then by definition of $r$ and $s$, the sets $\{N(h_i,s)\}_{i=1}^\infty, N(g,r)$ are mutually disjoint. 

 We give $A_\alpha$ coordinates $x,t$ where $t$ is the distance to the core geodesic $g$ and $x$  parametrizes the length about the core geodesic. Then
 \begin{eqnarray*}
\sum_{i}\int_{N(h_i,s)} e^{-2t} dA &\leq& \int_{A_\alpha\ssm N(g, r)} e^{-2t} dA \\
&= &2\int_0^{\ell_\alpha}\int_r^\infty e^{-2t} \cosh(t)dtdx\\ &=&  \ell_\alpha\left(e^{-r}+\frac{e^{-3r}}{3}\right).\end{eqnarray*}
To estimate the terms in the sum on the left we note that the integrals can be lifted to the hyperbolic plane and then by Lemma \ref{mean}
 $$e^{-2t_i} \leq \frac{1}{2\tan^{-1}(\sinh(s))\cosh^2(s) +2\sinh(s)}\int_{N(h_i,s)} e^{-2t} dA.$$
Therefore
$$\sum_{i}e^{-2t_i}  \leq \frac{\ell_\alpha\left(e^{-r}+\frac{e^{-3r}}{3}\right)}{2\tan^{-1}(\sinh(s))\cosh^2(s) +2\sinh(s)}.$$

Therefore by Riera's formula
$$\langle \grad \ell_\alpha,\grad \ell_\beta \rangle \leq \frac{2}{\pi}\ell_\alpha\left(\delta^\alpha_\beta +\frac{a(r+s)\left(e^{-r}+\frac{e^{-3r}}{3}\right)}{2\tan^{-1}(\sinh(s))\cosh^2(s) +2\sinh(s)}\right) = \frac{2}{\pi}\ell_\alpha\left(\delta^\alpha_\beta + G(r,s)\right)$$

As $G$ is the product of monotonically decreasing functions, it is monotonically decreasing. We now let  $\sinh(r) = 1/\sinh(\ell_\alpha/2)$ and $\sinh(s) = 1/\sinh(\ell_\beta/2)$ and define $F(\ell_\alpha,\ell_\beta) = G(r,s)$. 
Then 
$$e^{-r} = \frac{\sinh(\ell_\alpha/2)}{1+\cosh(\ell_\alpha/2)}$$
giving
\begin{eqnarray*}
F(\ell_\alpha,\ell_\beta) &=&
a(r+s)u(\ell_\alpha)v(\ell_\beta) \sinh(\ell_\alpha/2)\sinh^2(\ell_\beta/2)
\end{eqnarray*}
where 
  $$u(\ell_\alpha)= \frac{2\cosh(\ell_\alpha/2))+1}{3(\cosh(\ell_\alpha/2))+1)^2}\ \mbox{ and } \  v(\ell_\beta)= \frac{1}{\tan^{-1}\left(\csch(\ell_\beta/2)\right)\cosh^2(\ell_\beta/2)+\sinh(\ell_\beta/2)}.$$  
\eproof

We  now prove Theorem \ref{grad}.

{\bf Proof of Theorem \ref{grad}:}
We need to show that for $0 \leq z \leq w$ then 
$$F(z,w) \leq \frac{4}{3\pi}\sinh(z/2)\sinh^2(w/2).$$
We let $r = \sinh^{-1}(1/\sinh(z/2)), s=  \sinh^{-1}(1/\sinh(w/2))$. Then
$F(z,w) = a(r+s)u(z)v(w)\sinh(z/2)\sinh^2(w/2)$. We now show that $a(r+s)u(z)v(w) \leq 4/3\pi$ by showing it is maximized at $z = w = 0$.

 We first show $v$ is monotonically decreasing. We implicitly define $v_1(\sinh(w/2)) = 1/v(w)$. Then
 $$v_1(t) = (1+t^2)\tan^{-1}\left(\frac{1}{t}\right)+t$$
 giving
 $$v'_1(t) = 2t\tan^{-1}\left(\frac{1}{t}\right) + (1+t^2)\left(\frac{1}{1+\frac{1}{t^2}}\right)\left(\frac{-1}{t^2}\right)+1 = 2t\tan^{-1}\left(\frac{1}{t}\right).$$
  Therefore  $v_1$ is monotonically increasing, and $v$ is monotonically decreasing. It follows that $v(w) \leq v(0) = 2/\pi$.
 
 We now show $a(r+s)u(z) \leq 2/3$. By assumption $z\leq w$, giving $s \leq r$. Thus 
  $$u(z) =\frac{1}{4}\left(1+\frac{e^{-4r}}{3}\right)(1-e^{-2r}) \leq \frac{1}{4}\left(1+\frac{e^{-2(r+s)}}{3}\right)(1-e^{-2(r+s)}).$$ 
 We now use the expansion $\hat a (q) = \sum a_nq^{2n}$ from Lemma \ref{eRiera}. Letting $q = e^{-(r+s)}$ then $\hat a(q) = a(z)$ giving
 $$4a(r+s)u(z) \leq \hat a(q)\left(1+\frac{q^2}{3}\right)(1-q^2) = \left(\sum_{n=0}^\infty a_nq^{2n}\right)\left(1+\frac{q^2}{3}\right)(1-q^2) = \sum_{n=0}^\infty A_nq^{2n}.$$
 Computing we have
 $$A_n =(a_n-a_{n-1})+\frac{1}{3}(a_{n-1}-a_{n-2})$$
 where we define $a_{-1}=a_{-2} = 0$. For $n \geq 1$
 $$a_n-a_{n-1} = \left(\frac{2}{2n+3}+\frac{2}{2n+1}\right) -\left(\frac{2}{2n+1}+\frac{2}{2n-1}\right)= \frac{2}{2n+3}-\frac{2}{2n-1} < 0$$
  Thus $A_n < 0$ for $n \geq 2$. Also 
  $$A_1 = a_1-a_0+\frac{a_0}{3} =\frac{16}{15}-\frac{8}{3}+\frac{8}{9} < -\frac{32}{45}.$$
  It follows that $A_n < 0$ for all $n\neq 0$. Therefore $4a(r+s)u(z) \leq A_0 = 8/3$ giving
   $a(r+s)u(z)v(w) \leq 4/3\pi$. 
   \eproof

We define $F(t) =F(t,t)$. Then from above, we have the following;

\begin{corollary}
Let $S$ be a finite type hyperbolic surface and $\ell_\alpha$ be a geodesic length function for $\alpha$ simple. Then
$$\frac{2\ell_\alpha(X)}{\pi}\leq \|\nabla \ell_\alpha(X)\|^2 \leq \frac{2\ell_\alpha(X)}{\pi}\left(1+ F(\ell_\alpha(X))\right)$$
 where $F(t) \leq  (4/3\pi)\sinh^3(t/2).$
\label{lgrad}
\end{corollary}

\begin{comment}
{\bf Proof:}
We define  $F(t) = a(T)u(t)v(t)\sinh^{3}(t/2)$ where $\sinh(T/2)\sinh(t/2) = 1$. We need to calculate $a(T) = e^{2T}R(\cosh(T))$. We note that $\cosh(T/2) = \operatorname{cotanh}(t/2)$ giving
$$e^{T} = \frac{\cosh(t/2)+1}{\cosh(t/2)-1} \qquad \cosh(T) = \frac{\cosh^2(t/2)+1}{\sinh^2(t/2)}.$$
Thus
$$a(T) = \left(\frac{\cosh(t/2)+1}{\cosh(t/2)-1}\right)^2\left(\left(\frac{\cosh^2(t/2)+1}{\sinh^2(t/2)}\right)\log\left(\frac{\frac{\cosh^2(t/2)+1}{\sinh^2(t/2)}+1}{\frac{\cosh^2(t/2)+1}{\sinh^2(t/2)}-1}\right)-2\right).$$
Simplifying we get
$$a(T) = \frac{2(\cosh(t/2)+1)^2\left((1+\cosh^2(t/2))\log(\cosh(t/2))-\sinh^2(t/2)\right)}{(\cosh(t/2)-1)^2\sinh^2(t/2)}.$$
Combining with the formulae for $u,v$ gives the formula for $F$. As $u(0)= 1/4, v(0)=2/\pi$ and $\lim_{T\rightarrow \infty} a(T) = 8/3$. Expanding we have
$$F(t) = \frac{4}{3\pi}\sinh^3(t/2) + O(t^4) = \frac{t^3}{6\pi} + O(t^4).$$
Expanding $F$ at $t = 0$ further gives
$$F(t) =  \frac{t^3}{6\pi} -  \frac{t^5}{36\pi} + O(t^6).$$
 \eproof
\end{comment}

We note that Theorem \ref{grad} also gives a bound on $\|\nabla \ell_\alpha(X)\|$ in terms of collar radius. Defining $G(r) = G(r,r)$ then $G$ is monotonically decreasing with
$$G(r) = \frac{a(2r)\left(e^{-r}+\frac{e^{-3r}}{3}\right)}{2\tan^{-1}(\sinh(r))\cosh^2(r) +2\sinh(r)}.$$

\begin{corollary}
Let $S$ be a finite type hyperbolic surface and $\ell_\alpha$ be a geodesic length function for $\alpha$ simple. Let  $\alpha$ have an embedded neighborhood of radius $r_\alpha(X)$ in $X$. Then
$$\|\nabla \ell_\alpha(X)\|^2 \leq \frac{2\ell_\alpha(X)}{\pi}\left(1+ G(r_\alpha(X))\right).$$
Furthermore $G$ is monotonically decreasing with
$$F(t) = G\left(\sinh^{-1}\left(\frac{1}{\sinh(t/2)}\right)\right).$$
 \label{tgrad}\end{corollary}

From Corollary \ref{lgrad} the asymptotics of our bounds as $\ell_\alpha \to 0$ are easy to see. In particular, the difference between the upper and lower bounds is of order $\ell_\alpha^4$. In this form the asymptotics of our bounds are not as transparent when $\ell\to \infty$. For this purpose, it is useful to rephrase our bounds in terms of simpler functions.

Before doing so we first  state a theorem of Wolpert:
\begin{theorem}[{Wolpert, \cite{Wolpert:behavior}}]
Let $\ell_\alpha$ be a geodesic length functions on $\Teich(S)$, then
$$\|\grad \ell_\alpha(X)\|  \leq c\left(\ell_\alpha(X) + \ell_\alpha(X)^2e^{\frac{\ell_\alpha(X)}{2}}\right)$$
for some universal constant $c > 0$.
\label{wolpert_large}
\end{theorem}

Our bound gives an effective version of Wolpert's result with the same asymptotics as $\ell_\alpha\to \infty$.
\begin{corollary}
Let $\ell_\alpha$ be a geodesic length functions on $\Teich(S)$, then
$$\|\nabla \ell_\alpha (X)\|^2 \leq \frac{2}{\pi}\left(\ell_\alpha(X)+\frac{1}{3}\ell_\alpha(X)^2 e^{\ell_\alpha(X)/2}\right).$$
\end{corollary}

{\bf Proof:}
We have that the function  $F(t) = a(T)u(t)v(t)\sinh^{3}(t/2)$ where $T = 2\sinh^{-1}(1/\sinh(t/2))$.
Considering $u$ we have
  $$u(t)= \frac{2\cosh(t/2))+1}{3(\cosh(t/2))+1)^2} \leq \frac{2}{3\cosh(t/2)} \leq \frac{4e^{-t/2}}{3}.$$ 
 For $v(t)$ we consider $f(s) = \tan^{-1}(1/s) - 1/\sqrt{1+s^2}$ for $s > 0$. We have
 $$f'(s) = -\frac{1}{1+s^2} +\frac{s}{(1+s^2)^{3/2}} = \frac{s-\sqrt{1+s^2}}{(1+s^2)^{3/2}} \leq 0.$$
  Therefore $f$ is monotonically decreasing and $\lim_{s\rightarrow \infty} f(s) = 0$. It follows that for $s=\sinh(t/2)$ we get    $\tan^{-1}(\csch(t/2)) \geq \sech(t/2)$. Therefore
    $$v(t)= \frac{1}{\tan^{-1}\left(\csch(t/2)\right)\cosh^2(t/2)+\sinh(t/2)} \leq \frac{1}{\cosh(t/2)+\sinh(t/2)} = e^{-t/2}.$$  
We now bound $a(T)$. As $\sinh(T/2)\sinh(t/2) = 1$, we have 
$$e^{-T} = \frac{\cosh(t/2)-1}{\cosh(t/2)+1}.$$ 
By Lemma \ref{eRiera} we have the bound  $a(T) \leq 8/3-2\log(1-e^{-2T})$. Therefore
$$a(T) \leq \frac{8}{3}+2\log\left(\frac{(1+\cosh(t/2))^2}{4\cosh(t/2)}\right) = \frac{8}{3}+t + 2\log\left(\frac{(1+\cosh(t/2))^2}{4e^{t/2}\cosh(t/2)}\right).$$
As 
$$ \frac{(1+\cosh(t/2))^2}{4e^{t/2}\cosh(t/2)} = \frac{1}{4}\left(\frac{2}{e^{t}+1} + \frac{2}{e^{t/2}}+\frac{1+e^{-t}}{2}\right) \leq 1$$
we obtain $a(T) \leq 8/3 + t$.
Therefore
$$F(t) \leq  \frac{4}{3}\left(t+\frac{8}{3}\right)e^{-t}\sinh^3(t/2).$$ 
It follows that
$$\frac{F(t)}{te^{t/2}} \leq \frac{1}{6}\left(1+\frac{8}{3t}\right)(1-e^{-t})^3 \leq \frac{1}{6}+\frac{4}{9t}(1-e^{-t})^3= g(t).$$
By simple calculus, $g$ has a single critical point  $t_0 > 0$ that is the global maximum. Evaluating we get $g(t_0) \leq 1/3$. The result follows.
\eproof

\section{Bounding strata separation}
We now give an explicit bound on Wolpert's strata separation. Before doing so we prove the following elementary lemma.
\begin{lemma}
Let $M$ be a Riemannian manifold and $f\colon M \to \mathbb R$ be a smooth function.
Let $U$ and $L$ be non-negative integrable functions with
$$ L(f(x)) \le \|\nabla f (x)\| \le U(f(x))$$
for all $x\in M$.
Then if $x_t$ is an integral curve of $\nabla f$ that is defined on the interval $[a,b]$ we have
$$d(x_a, x_b) \le \int_{f(x_a)}^{f(x_b)} \frac1{L(s)}ds$$
and for any $x,y \in M$ with $f(x) \le f(y)$ we have
$$d(x,y) \ge \int_{f(x)}^{f(y)} \frac1{U(s)} ds.$$
\label{gradflow}
\end{lemma}

{\bf Proof:} We begin with the first inequality. We have
$$d\left(x_a,x_b \right) \leq  \operatorname{Length}\left(x_{[a,b]}\right) = \int_a^b \|\dot x_t\|dt = \int_a^b \| \grad f(x_t)\| dt.$$
If we make the substitution $s= f(x_t)$ we have
$$ds = df(\dot x_t)= \langle \grad f(x_t), \dot x_t\rangle dt =  \langle \grad f(x_t), \grad f(x_t)\rangle dt = \|\grad f(x_t)\|^2 dt$$
and therefore
$$d\left(x_a, x_b\right) =  \int_{f(x_a)}^{f(x_b)} \frac{1}{\| \grad f(x_t)\|} ds \leq  \int_{f(x_a)}^{f(x_b)} \frac{1}{L(f(x_t))} ds = \int_{f(x_a)}^{f(x_b)} \frac{1}{L(s)} ds.$$

Let $y_t$ be a smooth path in $M$ with $x=y_0$ and $y=y_1$.
Letting $s = f(y_t)$ we have
$$ds = df(\dot y_t) = \langle \grad f(y_t), \dot y_t\rangle dt \leq \|\dot y_t\|\cdot\|\grad f(y_t)\|dt.$$
We let $E \subset [0,1]$ where $s$ is monotonically increasing. Then
$$  \operatorname{Length}\left(y_{[0,1]}\right) = \int_0^1 \|\dot y_t\|dt \geq \int_E \frac{1}{\|\grad f(y_t)\|}dt  \geq \int_{f(x)}^{f(y)} \frac{1}{U(s)}ds.$$
As this holds for all paths from $x$ to $y$ we have
$$d(x,y) \ge \int_{f(x)}^{f(y)} \frac1{U(s)} ds.$$
\eproof

The following  proposition will allow us to apply this lemma to the gradient flow on length functions on $\Teich(S)$.
\begin{prop}\label{gradientflow}
Let $X_t$ be an integral curve of $\nabla \ell_\alpha$ and let $(a,b)$ be the maximal domain where $X_t$ is defined. Then
$$\underset{t\to a^+}{\lim }\  \ell_\alpha(X_t) = 0 \mbox{ and } \underset{t\to b^-}{\lim}\  \ell_\alpha(X_t) = +\infty.$$
Furthermore the limit of $X_t$ as $t\to a^+$ exists and lies in $\overline{\mathcal S_\alpha}$.
\end{prop}

\noindent
{\bf Proof:} By Theorem \ref{wgrad} an upper bound on $\ell_\alpha(X)$ gives an upper bound on $\|\nabla \ell_\alpha(X)\|$. Therefore if we fix $T\in (a,b)$ the length of the flow line $X_t$ on $(a,T]$ will be finite so $X_t$ converges to some $X_a \in \overline{\Teich(S)}$ as $t  \to a^+$. As $\nabla \ell_\alpha$ is non-zero on $\Teich(S)$ the limit must be in some boundary strata $\mathcal S_\tau$ where $\tau$ is a multi-curve on $S$. In particular if $\beta\subset \tau$ then $\underset{t\to a^+}{\lim}\ \ell_\beta(X_t) = 0$.

Note that for all $t\in (a,T]$ we have $\ell_\alpha(X_t) \le \ell_\alpha(X_T)$ so every curve on $X_t$ that intersects $\alpha$ will have length uniformly bounded away from zero by a constant depending on $\ell_\alpha(T)$. Therefore $\alpha$ and $\beta$ are disjoint if $\beta \subset \tau$.

We simplify notation and set $\ell_\alpha(t) = \ell_\alpha(X_t)$ and $\ell_\beta(t) = \ell_\beta(X_t)$. As $X_t$ is an integral curve of $\nabla \ell_\alpha$, $\ell'_\alpha(t)>0$. By the Riera formula (Theorem \ref{Riera}), the inner produce of $\nabla \ell_\alpha$ and $\nabla \ell_\beta$ is non-negative, so $\ell'_\beta(t) \ge 0$ and $\ell_\beta$ is non-decreasing. Therefore for $t\in(a,T]$ both $\ell_\alpha(t)$ and $\ell_\beta(t)$ are bounded above by $\max\{\ell_\alpha(T), \ell_\beta(T)\}$.
Again applying Theorem \ref{wgrad} we have
$$\ell'_\alpha(t) \ge \frac{2}\pi \ell_\alpha(t) \mbox{ and } \ell'_\beta(t) \le C\ell_\alpha(t)^2 \ell_\beta(t)^2$$
for $t \in (a,T]$ where $C$ depends on $\max\{\ell_\alpha(T), \ell_\beta(T)\}$.

If $\ell_\beta(X_a) = 0$ and $\ell_\alpha(X_a) = \epsilon>0$
 then
$$\underset{t\to a^+}{\lim}\ \log \frac{\ell_\beta(t)}{\ell_\alpha(t)} \to -\infty.$$
However, when $\ell_\beta(t) < \frac2{C\epsilon^2\pi}$ this function is decreasing as
\begin{eqnarray*}
\left(\log \frac{\ell_\beta(t)}{\ell_\alpha(t)}\right)' &=& \frac{\ell'_\beta(t)}{\ell_\beta(t)} - \frac{\ell'_\alpha(t)}{\ell_\alpha(t)} \\
& \le & C \ell_\alpha(t)^2\ell_\beta(t) - \frac2\pi\\
& <& 0.
\end{eqnarray*}
This contradiction proves that $\ell_\alpha(X_a) = 0$ and therefore $X_a \in \overline{\mathcal S_\alpha}$.

For the second statement, if the limit of $\ell_\alpha(X_t)$ as $t\to b^+$ is finite
 then, as above, the integral curve will have finite length and must have a limit in some boundary strata $\mathcal S_\tau \subset \overline{\Teich(S)}$. However, if $\tau$ intersects $\alpha$ then the length of $\alpha$ will be infinite in the limit, a contradiction. Therefore $\tau$ must be disjoint from $\alpha$. However, by the Riera formula, the length of every curve disjoint from $\alpha$ will increase along $X_t$, again a contradiction. This establishes the second claim.
\eproof

In the following, as the surface $S$ is understood, we will  denote strata as $\mathcal S_\tau$ where $\tau$ is a multicurve.

Motivated by Theorem \ref{grad} we define 
$$H(a,b) = \int_a^b \frac{dt}{\sqrt{\frac{2t}{\pi}\left(1+F(t)\right)}} \qquad \mbox{ and }  \qquad K(a,b) = \int_a^b \frac{dt}{\sqrt{\frac{2t}{\pi}}} = \sqrt{2\pi b} - \sqrt{2\pi a}.$$
We will often be interested in the case when $a=0$ and in this case we will write $H(b) = H(0,b)$ and $K(b) = K(0,b)$.
We denote  the level sets of the length function $\ell_\alpha$ by 
$$\mathcal S_\alpha^L = \ell_\alpha^{-1}(L) \subseteq\overline{\Teich(S)}.$$
Combining Lemma \ref{gradflow} and Proposition \ref{gradientflow} to the bounds in Theorem \ref{grad} we get:
\begin{theorem}
Let $\alpha$ be a simple closed curve on $S$. Then for $a,b\in [0,\infty)$ and $X\in {\mathcal S}_\alpha^a$ we have 
$$|H(a,b)| \leq d_{\rm WP}(\mathcal S_\alpha^a, \mathcal S_\alpha^b)\le d_{\rm WP}(X, \mathcal S_\alpha^b) \leq |K(a,b)|.$$
\label{levelcurve}\end{theorem}

\subsubsection*{Accuracy of  bounds}
We now discuss the accuracy of our bounds. For this purpose we define functions
$$\mathcal D^+(\ell) = \underset{\{X \mbox{ with } \ell_\alpha(X) = \ell\}}{\sup} \| \nabla \ell_\alpha(X)\|^2$$
and
$$\mathcal D^-(\ell) = \underset{\{X \mbox{ with } \ell_\alpha(X)=\ell\}}{\inf} \| \nabla \ell_\alpha(X)\|^2.$$
It is not hard to check that $\mathcal D^-(\ell)  = \frac{2\ell}\pi$ and therefore the lower bound is optimal. In particular, one can find a sequence $X_i$ where $\ell_\alpha(X_i) = \ell$ and the widths of the maximal collars about $\alpha$ on $X_i$ go to infinity. By Corollary \ref{tgrad} as the width limits to infinity the difference between the upper and lower bound will converge to zero.

To estimate $\mathcal D^+(\ell)$ we return to the family of rectangular punctured tori from the proof Lemma \ref{dbound}. Here there are two curves $\alpha$ and $\beta$ meeting orthogonally with $\sinh(\ell_\alpha/2)\sinh(\ell_\beta/2) =1$. Then
$$|\grad \ell_\alpha|^2 \geq \frac{2}{\pi}\left(\ell_\alpha +R(\cosh(\ell_\beta))\right) \geq \frac{2}{\pi}\left(\ell_\alpha +\frac{2}{3}\left(\frac{1}{\cosh^2(\ell_\beta)}\right)\right).$$
Thus
$$|\grad \ell_\alpha|^2  \geq \frac{2}{\pi}\left(\ell_\alpha +\frac{2}{3}\frac{\sinh^4(\ell_\alpha/2)}{(1+\cosh^2(\ell_\alpha/2))^2} \right).$$
We consider $\ell_\alpha$ small. Then
$$|\grad \ell_\alpha|^2  \geq \frac{2}{\pi}\left(\ell_\alpha +\frac{\ell_\alpha^4}{24} +O(\ell_\alpha^6)\right).$$
We note that by Corollary \ref{lgrad} the upper bound for $\ell_\alpha$ small gives
$$|\grad \ell_\alpha|^2 \leq \frac{2\ell_\alpha}{\pi}\left(1+ F(\ell_\alpha)\right) = \frac{2}{\pi}\left(\ell_\alpha +\frac{\ell_\alpha^4}{6\pi} +O(\ell_\alpha^6) \right).$$
Thus for short geodesics $\mathcal D^+(\ell)$  and our upper bound differ at order 4.

Similarly we consider $\ell_\alpha$ large. As $\sinh(\ell_\alpha/2)\sinh(\ell_\beta/2) =1$ differentiating we have
$$\|\grad \ell_\alpha\|^2 = \sinh^2(\ell_\alpha/2)\|\grad \ell_\beta\|^2 \geq \sinh^2(\ell_\alpha/2)\frac{2}{\pi}\left(\ell_\beta +\frac{\ell_\beta^4}{24} +O(\ell_\beta^6)\right)$$
As $\ell_\alpha$ is large 
$$\sinh(\ell_\alpha/2) = \frac{e^{\ell_\alpha/2}}{2}(1+ O(e^{-\ell_\alpha}))\qquad \qquad \sinh(\ell_\beta/2) = \frac{\ell_\beta}{2}(1+ O(\ell_\beta^2)) = \frac{\ell_\beta}{2}(1+ O(e^{-\ell_\alpha})).$$
 As $\sinh(\ell_\alpha/2)\sinh(\ell_\beta/2)=1$
then $\ell_\beta = 4e^{-\ell_\alpha/2}(1+ O(e^{-\ell_\alpha}))$ giving
$$\|\grad \ell_\alpha\|^2 \geq \sinh^2(\ell_\alpha/2)\frac{2}{\pi}\left(\ell_\beta +\frac{\ell_\beta^4}{24} +O(\ell_\beta^6)\right) = \frac{2}{\pi}e^{\ell_\alpha/2}\left(1 + O(e^{-\ell_\alpha})\right).$$
We note that the upper bound is
$$\|\grad \ell_\alpha\|^2 \leq \frac{1}{3\pi}\ell_\alpha^2e^{\ell_\alpha/2}\left(1 + O(1/\ell_\alpha)\right).$$
Thus as $\ell$ goes to infinity, $\mathcal D^+(\ell)$ grows of order at least $e^{\ell/2}$ while our upper bound grows of order $e^{\ell/2+\epsilon}$.

\section{Orthogonal projection onto strata}
\newcommand{\cat}{{\operatorname{CAT}}}
The Weil-Petersson completion $\overline{\Teich(S)}$ is a $\cat(0)$ space. 
Let $\tau$ be a multicurve in $S$, $\mathcal S_\tau$ the associated strata and $S_\tau = S\smallsetminus \tau$. Then $\mathcal S_\tau$ is isometric to $\Teich(S_\tau)$ and the closure $\overline{\mathcal S_\tau}$ is convex in $\overline{\Teich(S)}$ (see \cite{Yamada:cat}, \cite{Wolpert:Nielsen}). Note that if $S_\tau$ is disconnected then $\Teich(S_\tau)$ is the product of the Teichm\"uller spaces of each component.

Now, let $\tau_0$ and $\tau_1$ be multicurves in $S$ and $\mathcal S_{\tau_0}$ and $\mathcal S_{\tau_1}$ the associated strata. We will show that the infimum of distance between $\mathcal S_{\tau_0}$ and $\mathcal S_{\tau_1}$ is attained on any stratum $\mathcal S_\sigma$ for which is $\sigma$ is mutually disjoint from both $\tau_0$ and $\tau_1$. Specifically we prove:

\begin{theorem}
Let $\tau_0, \tau_1$, and $\sigma$ be multicurves with $i(\tau_i, \sigma) = 0$ for $i=0,1$. If $\hat \tau_i = \tau_i \cup \sigma$ then
$$d_{\rm WP}(\mathcal S_{\tau_0},\mathcal S_{\tau_1}) = d_{\rm WP}
(\mathcal S_{\hat\tau_0},\mathcal S_{\hat\tau_1}).$$
\label{optimal}
\end{theorem}

In a $\cat(0)$ space the nearest point projection to a convex set is $1$-Lipschitz (see \cite[Proposition 2.4]{Bridson:Haefliger:npc}). 
Here we will project to the closure $\overline{\mathcal S_\sigma}$ and the theorem will follow once we show that this projection maps $\mathcal S_{\tau_i}$ into $\overline{\mathcal S_{\hat\tau_i}} \subset \overline{\mathcal S_{\tau_i}}$. This in turn follows quickly from Wolpert's characterization of tangent cones in the Weil-Petersson metric (see \cite{Wolpert:behavior}). We begin by reviewing this work.

Given $p, q, r \in \overline{\Teich(S)}$ we let $\angle(p;q,r)$ be the angle at $p$ in the comparison Euclidean triangle with side lengths $d_{\rm WP}(p,q)$, $d_{\rm WP}(q,r)$ and $d_{\rm WP}(p,r)$. Let $b(t)$ and $c(t)$ be constant speed geodesic segments starting at $p$. The $\cat(0)$ property implies that if $0< s_0 \le s_1$ and $0<t_0\le t_1$ then
$$\angle(p; b(s_0), c(t_0)) \le \angle(p; b(s_1), c(t_1))$$
and therefore
$$\angle(b, c) = \underset{t\to 0}{\lim}\angle(p; b(t), c(t))$$
is defined. Let $|b|$ and $|c|$ be the (constant) speed of the two segments. We define an equivalence relation where $b\sim c$ if $|b| = |c|$ and $\angle(b,c) = 0$. If we take all geodesic segments beginning at $p$ and take the quotient under this equivalence relation we have the {\em Alexandrov tangent cone} at $p$. At points in $\Teich(S)$ this is the usual tangent space at $p$. 

We also define an {\em inner product} by
$$\langle b,c\rangle = |b|\cdot |c| \cos(\angle(b,c)).$$
\begin{theorem}[Wolpert, \cite{Wolpert:behavior}]\label{cone}
Let $\tau = \{\gamma_1, \dots, \gamma_k\}$ be a multicurve and assume that $p \in \mathcal S_\tau$. The the Alexandrov tangent cone at $p$ is
$$\mathbb R^{|\tau|}_{\ge 0} \times T_p\mathcal S_\tau$$
where the inner product is the product of the standard inner produce on $\mathbb R^{|\tau|}$ and the Weil-Petersson inner product on $T_p\mathcal S_\tau$.
Furthermore if $b(t)$ is a constant speed geodesic segment starting at $p$ and $\ell_{\gamma_i}(b(t)) = 0$ then the $i$th coordinate of $b$ in the tangent cone is zero.
\end{theorem}

Given a multicurve $\sigma$ let
$$\pi_\sigma\colon \overline{\Teich(S)} \to \overline{\mathcal S_\sigma}$$
be the nearest point projection.
\begin{lemma}
Let $\sigma$ be a multicurve in $S$ and $p$ and $q$ points in $\overline{\Teich(S)}$ with $p = \pi_\sigma(q)$. Then $p \in \mathcal S_{\hat\sigma}$ where $\hat\sigma$ is a (possibly trivial) extension of $\sigma$.
Let $b(t)$ be a geodesic segment from $p$ to $q$. Then the image of $b$ in the tangent cone is orthogonal to ${\mathbb R}^{|\hat\sigma\smallsetminus\sigma|}_{\ge 0}\times T_p\mathcal S_{\hat\sigma}$.
\end{lemma}

{\bf Proof:} Let $c\colon (-\epsilon, \epsilon) \to \mathcal S_{\hat\sigma} \subset \overline{\mathcal S_\sigma}$ be a constant speed geodesic with $c(0) = p$. If we let $\bar c(t) = c(-t)$ then $\angle(c, \bar c) = \pi$. By (3) of \cite[Proposition 2.4]{Bridson:Haefliger:npc} the angles $\angle(b, c)$ and $\angle(b, \bar c)$ are at least $\pi/2$. Therefore they must be equal to $\pi/2$ and hence $b$ is orthogonal to $T_p\mathcal S_{\hat\sigma}$. In particular, by Theorem \ref{cone}, $b$ lies in $\mathbb R^{|\hat\sigma|}_{\ge 0}$.

Every vector in ${\mathbb R}^{|\hat\sigma\smallsetminus \sigma|}_{\ge 0}$ is represented by a geodesic segment $c\colon [0,\epsilon) \to \overline{\mathcal S_\sigma}$ with $c(0) = p$. In particular $d_{\rm WP}(q, c(t)) > d_{\rm WP}(q, p)$ for all $t\in (0,\epsilon)$. As above, (3) of \cite[Proposition 2.4]{Bridson:Haefliger:npc} implies that $\angle(b,c) \ge \pi/2$. However, as $b$ lies in ${\mathbb R}^{|\hat\sigma|}_{\ge 0}$, we must have that $\angle(b,c) = \pi/2$.  \eproof

\begin{prop}\label{projection}
Let $\tau$ and $\sigma$ be multicurves with $i(\tau, \sigma)= 0$ and let $\hat\tau = \tau \cup \sigma$. Then
$$\pi_\sigma(\mathcal S_\tau) \subset \overline{\mathcal S_{\hat\tau}}.$$
\end{prop}

{\bf Proof:} Let $q$ be a point in $\mathcal S_\tau$ and $p = \pi_{\hat\tau}(q)$ and $r = \pi_\sigma(q)$ its nearest point projections to $\overline{\mathcal S_{\hat\tau}}$ and $\overline{\mathcal S_\sigma}$. By the previous lemma the angles of the triangle $qpr$ at $p$ and $r$ are $\pi/2$ so in the Euclidean comparison triangles the corresponding angles must be at least $\pi/2$. However, if $p\neq r$ then the angle at $q$ in the comparison triangle will be $>0$, a contradiction.\eproof
\\
{\bf Proof of Theorem \ref{optimal}:} As $\mathcal S_{\hat{\tau}_i}$ is contained in $\overline{\mathcal S_{\tau_0}}$ we have
$$d_{\rm WP}(\mathcal S_{\tau_0},\mathcal S_{\tau_1}) \le d_{\rm WP}
(\mathcal S_{\hat\tau_0},\mathcal S_{\hat\tau_1}).$$
On the other hand, for any $X_0\in \mathcal S_{\tau_0}$ and $X_1 \in \mathcal S_{\tau_1}$ we have 
$$d_{\rm WP}(X_0, X_1) \ge d_{\rm WP}(\pi_\sigma(X_0), \pi_\sigma(X_1))$$
as the nearest point projection is $1$-Lipschitz. By Proposition \ref{projection}, $\pi_\sigma(X_i) \subset \overline{\mathcal S_{\hat{\tau}_i}}$ so
$$d_{\rm WP}(\mathcal S_{\tau_0},\mathcal S_{\tau_1}) \ge d_{\rm WP}
(\overline{\mathcal S_{\hat\tau_0}}, \overline{\mathcal S_{\hat\tau_1}})=d_{\rm WP}
(\mathcal S_{\hat\tau_0},\mathcal S_{\hat\tau_1}).$$
\eproof

\section{Topological properties of nearby strata}

We now prove Theorem \ref{smallsep} which we first restate. \newline

\noindent{\bf Theorem \ref{smallsep}}
{\em Let  $\mathcal S_\sigma,\mathcal S_\tau$ be two strata in $\Teich(S)$. Then one of the following holds;
  \begin{enumerate}
\item  $i(\sigma,\tau) = 0$ and  $d_{\rm WP}(\mathcal S_\sigma,\mathcal S_\tau)  = 0.$
\item $i(\sigma,\tau) = 1$ and  $d_{\rm WP}(\mathcal S_\sigma,\mathcal S_\tau)  = \delta_{1,1}.$
\item $i(\sigma,\tau) > 1$ and  $d_{\rm WP}(\mathcal S_\sigma,\mathcal S_\tau) \geq 7.61138.$ 
 \end{enumerate}
}

{\bf Proof:}
If  $i(\sigma,\tau) = 0$ then the closures of the strata intersect and therefore $d_{\rm WP}(\mathcal S_\sigma,\mathcal S_\tau) = 0$. 

Now assume that $i(\sigma, \tau)=k >0$ and that for every $\alpha \in \sigma$ we have $i(\alpha, \tau) = 0$ or $1$. Note that this implies that for every $\beta \in \tau$ then $i(\beta, \sigma) =0$ or $1$ and if $i(\sigma,\tau) = 1$ this condition automatically holds. Then the surface filled by $\sigma$ and $\tau$ will be a collection of punctured tori and annuli. Let $\mu$ be a maximal multicurve such that $i(\sigma, \mu) = i(\tau, \mu)=0$. Then $S\smallsetminus\mu$ will be the collection of $k$ punctured tori filled by $\sigma$ and $\tau$ along with a collection of thrice punctured spheres. If we let $\hat\sigma = \sigma\cup \mu$ and $\hat\tau = \tau\cup \mu$ then by Theorem \ref{optimal}
$$d_{WP}(\mathcal S_\sigma, \mathcal S_\tau) = d_{WP}(\mathcal S_{\hat\sigma}, \mathcal S_{\hat\tau}).$$
The strata $\mathcal S_{\hat\sigma}$ and $\mathcal S_{\hat\tau}$ are both maximal and hence each are a single point. As $\mu$ is a multicurve contained in both $\hat\sigma$ and $\hat\tau$, these strata are in the closure of $\overline{\mathcal S_\mu}$.  Furthermore $\overline{\mathcal S_\mu}$ is the product of $k$ copies of the Weil-Petersson completion of the Teichm\"uller space of the punctured torus and when we project to each factor the image of the strata $\mathcal S_{\hat\sigma}$ and $\mathcal S_{\hat\tau}$ are curves intersecting once. It follows that
$$d_{WP}(\mathcal S_{\hat\sigma}, \mathcal S_{\hat\tau}) = \sqrt k\delta_{1,1}.$$
Therefore if  $i(\sigma, \tau) = 1$ we have
$$d_{\rm WP}(\mathcal S_\sigma,\mathcal S_\tau) = \delta_{1,1}$$
and if $i(\sigma, \tau) = k \ge 2$ then by Lemma \ref{dbound}
$$d_{\rm WP}(\mathcal S_\sigma,\mathcal S_\tau) \geq \sqrt{2}\delta_{1,1} > \lbs.$$

Now we can assume, without loss of generality, that there is a curve $\alpha \in \sigma$ and curves $\beta_1$ and $\beta_2$ in $\tau$ (possibly with $\beta_1 =\beta_2$) and $i(\alpha, \beta_1\cup \beta_2) \ge 2$. Let $c$ be any path from $\mathcal S_\sigma$ to $\mathcal S_\tau$ and choose $t_0$ such that at $c(t_0) =X$ we have $\max\{\ell_{\beta_1}(X), \ell_{\beta_2}(X)\} = 2\epsilon_2$ where $\epsilon_2$ is the Margulis constant in dimension two. Therefore the collars about $\beta_1$ and $\beta_2$ have length at least $2\epsilon_2$ and as $i(\alpha,\beta_1\cup\beta_2) = 2$ this implies  $\ell_\alpha(X) \geq 4\epsilon_2$. Then by Theorem \ref{levelcurve}, $d_{\rm WP}(X,\mathcal S_\sigma) \geq  H(4\epsilon_2) $ and $d_{\rm WP}(X, \mathcal S_\tau) \geq  H(2\epsilon_2).$  Thus 
 $$d_{\rm WP}(\mathcal S_\sigma,\mathcal S_\tau) \geq H(4\epsilon_2) +H(2\epsilon_2).$$
 Evaluating we obtain $H(4\epsilon_2) +H(2\epsilon_2) \geq 7.61138.$
Thus if $i(\sigma,\tau) > 1$ and $d_{\rm WP}(\mathcal S_\sigma,\mathcal S_\tau) \geq 7.61138$.
 \eproof

\subsubsection*{Topology of supporting surface}

If the subsurface $S(\sigma, \tau)\subset S$ filled by $\sigma$ and $\tau$ has $n>1$ non-annular components then by the above 
$$d_{\rm WP}(\mathcal S_\sigma,\mathcal S_\tau)  \geq \sqrt{2} \delta_{1,1} > \lbs.$$
Thus if $d_{\rm WP}(\mathcal S_\sigma,\mathcal S_\tau) \leq \lbs
 $ then $S(\mu,\tau)$ has a  single non-annular component.
Also by the above, if $d_{\rm WP}(\mathcal S_\sigma,\mathcal S_\tau) \leq 7.61138$ then the non-annular component is a punctured torus with $i(\mu,\tau) = 1$ and  in fact $d_{\rm WP}(\mathcal S_\sigma,\mathcal S_\tau) = \delta_{1,1}$.

\subsubsection*{Separating curves and punctured spheres \label{nps}}
The above shows that for any finite type surface, $\delta_{1,1}$ is a lower bound on the  distance between strata in $\Teich(S)$ whose closures do not intersect. Also it follows that it is attained for any $S$ with a non-separating curve. The only case left is the  $n$-punctured sphere $S_{0,n}$ for $n \geq 4$. For completeness, we now consider this case.

In a punctured sphere every curve is separating so any two curves with non-trivial intersection will intersect an even number of times. In particular, on the 4-punctured sphere any two distinct curves intersect and the minimal intersection is two. In parallel with the punctured torus case, if $\alpha$ and $\beta$ are simple closed curves in $S_{0,4}$ with $i(\alpha,\beta) = 2$ we define
$$\delta_{0,4} = d_{\rm WP}(\mathcal S_\alpha(S_{0,4}),\mathcal S_\beta(S_{0,4})).$$

We note that there is an canonical isomorphism between $\overline{\Teich(S_{1,1})}$ and $\overline{\Teich(S_{0,4})}$ and as the area of 4-punctured hyperbolic spheres is twice that of punctured tori this isomorphism scales the Weil-Petersson metric by the $\sqrt2$. Two noded surfaces in $\overline{\Teich(S_{1,1})}$ whose nodes intersect once will be taken to noded surfaces in $\overline{\Teich(S_{0,4})}$ where the nodes intersect twice and therefore
$$\delta_{0,4} = \sqrt2\delta_{1,1}$$
Therefore by the bounds on $\delta_{1,1}$ in Lemma \ref{dbound} we have $\delta_{0,4} \in (\lbs, \ubs).$

The usual collar lemma states that if $\alpha$  is a simple closed geodesic in a complete hyperbolic surface $X$ then $\alpha$ has an embedded collar of width $r$ with $\sinh(r/2) = 1/\sinh(\ell_\alpha(X))$. If $\alpha$ is non-separating then this result is optimal: for any $\epsilon>0$ there is a hyperbolic structure $X$ (on any hyperbolizable surface $S$) such that $\alpha$ doesn't have a collar of width $r+\epsilon$. However, for separating curves this can be improved. While the proof is elementary we were unable to find a reference so we include one here. (See \cite{parlier:collar} for a similar observation.)
\begin{lemma}
Let $\alpha$ be a separating curve on a complete hyperbolic surface $X$. Then $\alpha$ has an embedded collar of width $r$ with

$$\sinh(\ell_\alpha(X)/4) \sinh(r/2) \geq 1.$$
\label{collar2}
\end{lemma}

\begin{figure}[htbp] 
   \centering
   \includegraphics[width=4in]{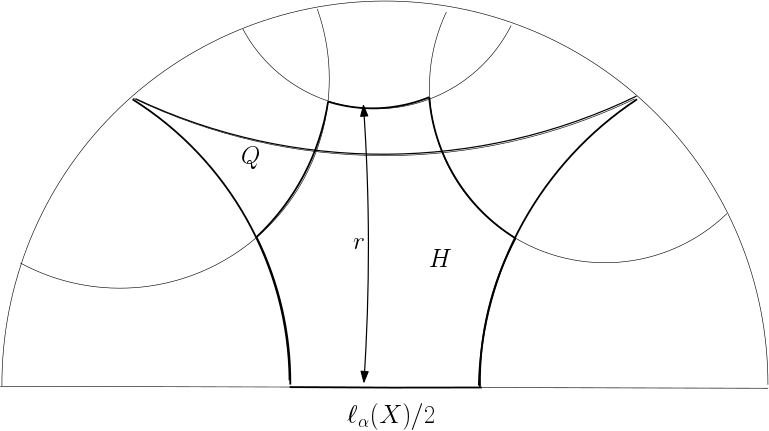} 
   \caption{$r > 2\sinh^{-1}(1/\sinh(\ell_\alpha(X)/4))$}
   \label{hex}
\end{figure}

{\bf Proof:} Let $\beta$ be the shortest non-trivial geodesic arc from $\alpha$ to itself. Then we can choose $r$ to be  the length of $\beta$. As $\alpha$ is separating, $\beta$ starts and ends on the same side of $\alpha$. 
 Therefore $\alpha$ and $\beta$ are  supported  on a pair of pants $P$ in $X$. We decompose $P$ into two isometric right-angled hexagons in the standard way by taking perpendiculars between boundary components of $P$.
This hexagon has base of length $\ell_\alpha(X)/2$.  We extend the sides of $H$ to geodesics in $\Hp$. The sides perpendicular to the base are distance $\ell_\alpha(X)/2$ apart and therefore are the opposite sides of an ideal quadrilateral $Q$ with the two other sides  a distance $2\sinh^{-1}(\ell_\alpha(X)/4)$ apart (see Figure \ref{hex}).
 The  geodesic opposite the base geodesic is separated from the base geodesic by a side of $Q$. Therefore the distance from the base to the opposite geodesic is at least $\sinh^{-1}(1/\sinh(\ell_\alpha(X)/4))$.

As $\beta$ is the union of two geodesic arcs joining the base of $H$ to its opposite side and $r$ is  the length of $\beta$, we have
$$r \geq 2\sinh^{-1}\left(\frac{1}{\sinh(\ell_\alpha(X)/4)}\right).$$
\eproof

In the usual collar lemma, the standard collars are disjoint. We emphasize that this does not hold for the collars we construct here.

Using the above we can improve our gradient bound for separating curves. We have
\begin{theorem}
Let $S$ be a finite type surface and $\ell_\alpha$ be a geodesic length function for $\alpha$ a simple separating curve on $S$. Then for $X \in \Teich(S)$
$$\|\nabla \ell_\alpha(X)\|^2 \leq \frac{2\ell_\alpha(X)}{\pi}\left(1+ F(\ell_\alpha(X)/2)\right).$$
Furthermore 
$$d_{\rm WP}(\mathcal S^a_\alpha, \mathcal S^b_\alpha) \geq H_s(a,b)$$
where
$$H_s(a,b) = \int_a^b \frac{dt}{\sqrt{\frac{2t}{\pi}\left(1+ F(t/2)\right)}}.$$
\label{sgrad}
\end{theorem}

{\bf Proof:}
The proof is the same as in Theorem \ref{grad}. The only difference is that the embedded neighborhood 
has width $2\sinh^{-1}(1/\sinh(\ell_\alpha(X)/4))$ rather than $2\sinh^{-1}(1/\sinh(\ell_\alpha(X)/2))$. Thus we can substitute $\ell_\alpha(X)/2$ into the lower  bound in Corollary \ref{lgrad} to obtain the new lower bound. We note the linear factor arises from integrating in the $\alpha$ direction in the collar and therefore remains unchanged. The  Weil-Petersson distance bound   follows immediately as in Lemma \ref{levelcurve}.
\eproof

We repeat the proof of Theorem \ref{smallsep} for the punctured sphere case. For simplicity, we will let $H_s(t) = H_s(0,t)$.

 \begin{theorem}
  Let  $\mathcal S_\sigma(S),\mathcal S_\tau(S)$ be two strata in $\Teich(S)$ for $S$ an $n$-punctured sphere. Then one of the following holds;
  \begin{enumerate}
\item  $i(\sigma,\tau) = 0$ and  $d_{\rm WP}(\mathcal S_\sigma(S),\mathcal S_\tau(S))  = 0.$
\item $i(\sigma,\tau) = 2$ and  $d_{\rm WP}(\mathcal S_\sigma(S),\mathcal S_\tau(S))  = \delta_{0,4}.$
\item $i(\sigma,\tau) > 2$ and  $d_{\rm WP}(\mathcal S_\sigma(S),\mathcal S_\tau(S))> 10.09656.$ 
 \end{enumerate}
 \label{smallsep2}
\end{theorem}

{\bf Proof:}
If  $i(\sigma,\tau) = 0$ then the closures of the strata intersect and therefore $d_{\rm WP}(\mathcal S_\sigma,\mathcal S_\tau) = 0$. 

Now assume that $i(\sigma, \tau)=2k >0$ and that for every $\alpha \in \sigma$ we have $i(\alpha, \tau) = 0$ or $2$. Then by the same argument as in Theorem \ref{smallsep} we can decompose into 4-punctured spheres and get
$$d_{WP}(\mathcal S_{\sigma}, \mathcal S_{\tau}) = \sqrt k\delta_{0,4}.$$
Therefore if  $i(\sigma, \tau) = 2$ we have
$$d_{\rm WP}(\mathcal S_\sigma,\mathcal S_\tau) = \delta_{0,4}$$
and if $i(\sigma, \tau) = 2k \ge 4$ then by Lemma \ref{dbound}
$$d_{\rm WP}(\mathcal S_\sigma,\mathcal S_\tau) \geq \sqrt{2}\delta_{0,4}   =2\delta_{1,1} > 13.145 .$$

Now we can assume one of the following;
\begin{itemize}
\item there is curve $\alpha \in \sigma$ and curve $\beta \in \tau$ with $i(\alpha, \beta) \ge 4$.
\item there is curve $\alpha \in \sigma$ and curves $\beta_1, \beta_2 \in \tau$ and $i(\alpha, \beta_1)= i(\alpha, \beta_2) = 2$.
\end{itemize}

 In the first case, we let $c$ be any path from $\mathcal S_\sigma$ to $\mathcal S_\tau$ and choose $t_0$ such that at $c(t_0) =X$ we have $\ell_\beta(X) = L$. Therefore by Lemma \ref{collar2} above, $\alpha$ has an embedded collar of width $2\sinh^{-1}(1/\sinh(L/4))$. Therefore   $\ell_\alpha(X) \geq 8\sinh^{-1}(1/\sinh(L/4))$. Then by Theorem \ref{sgrad}, 
 $$d_{\rm WP}(\mathcal S_\sigma,\mathcal S_\tau) \geq H_s(L) +H_s\left(8\sinh^{-1}\left(\frac{1}{\sinh(L/4)}\right)\right) = W_1(L)$$
We  choose $L=3.678$ and
  evaluating we get
 $$d_{\rm WP}(\mathcal S_\sigma,\mathcal S_\tau) \geq W_1(3.678) \geq 10.76596.$$

In the second case,  we choose $t_0$ such that at $c(t_0) =X$  and
 $L = \max\{\ell_{\beta_1}(X),\ell_{\beta_2}(X)\}$. 
Then $\beta_1\cup\beta_2$ split $\alpha$ into 4 geodesic arcs with endpoints in $\beta_1\cup\beta_2$.
Two of the arcs have endpoints in the same component of $\beta_1\cup\beta_2$ and  therefore by Lemma \ref{collar2} are both of length at least $2\sinh^{-1}(1/\sinh(L/4))$.  The other two geodesic arcs have one endpoint in $\beta_1$ and another in $\beta_2$. Then using the fact that the  collars about $\beta_1,\beta_2$ of width $2\sinh^{-1}(1/\sinh(L/2))$ are disjoint we have each of these arcs are of length at least $2\sinh^{-1}(1/\sinh(L/2))$. Thus
$$\ell_\alpha(X) \geq 4\sinh^{-1}\left(\frac{1}{\sinh(L/4)}\right)+4\sinh^{-1}\left(\frac{1}{\sinh(L/2)}\right).$$
Thus
$$d_{\rm WP}(\mathcal S_\sigma,\mathcal S_\tau) \geq H_s(L) + H_s\left(4\sinh^{-1}\left(\frac{1}{\sinh(L/4)}\right)+4\sinh^{-1}\left(\frac{1}{\sinh(L/2)}\right)\right) = W_2(L).$$
We choose $L = 2.420$ and get
$$d_{\rm WP}(\mathcal S_\sigma,\mathcal S_\tau) \geq W_2(2.42) \geq 10.09656.$$
Thus if $i(\sigma,\tau) > 2$ then $d_{\rm WP}(\mathcal S_\sigma,\mathcal S_\tau) \geq 10.09656$.
\eproof

\subsubsection*{Strata distances and gaps}
From the above, if $S$ has positive genus then the minimal distance between strata $\mathcal S_\sigma,\mathcal S_\tau$ with $i(\sigma,\tau) \neq 0$ is $\delta_{1,1}$ and is achieved if and only if $i(\sigma,\tau) = 1$. 
Furthermore  if $i(\sigma,\tau) > 1$   then the distance between the strata is at least $H(4\epsilon_2)+H(2\epsilon_2) $. Therefore there is a gap in the distances from $\delta_{1,1}$ to $H(4\epsilon_2)+H(2\epsilon_2) $ of size
$$H(4\epsilon_2)+H(2\epsilon_2) - \delta_{1,1} \geq 7.61138- \ub =  0.95535.$$ 

Similarly if $S$ is an n-punctured sphere with $n \geq 4$,  then the minimal distance between strata $\mathcal S_\sigma,\mathcal S_\tau$ with $i(\sigma,\tau) \neq 0$ is $\delta_{0,4}$ and is achieved if and only if $i(\sigma,\tau) = 2$. Furthermore  if $i(\sigma,\tau) > 2$   then the distance between the strata is at least $W_2(2.42)$. Therefore there is a gap in the distances from $\delta_{0,4}$ to $W_2(2.42)$ of size

$$W_2(2.420) - \delta_{0,4} \geq 10.09656-\ubs  =  .68351.$$

\section{Gradient bounds at systoles and the in-radius of $\Teich(S)$}
A {\em systole} is a shortest closed geodesic on a Riemannian manifold.
The {\em systole function} 
$$\ell_{\sys}: \Teich(S) \rightarrow \R_{> 0}$$
is the length of the systole at $X \in \Teich(S)$. The systole function is a proper, bounded function to $(0,\infty)$ (as it extends continuously to zero on $\partial\overline{\T(S)}$) and therefore
$$\sys(S) = \max_{X\in \Teich(S)}\ell_{\rm sys}(X)$$
is defined. Note that for a fixed curve $\alpha$ we have bounded from below the distance between $X$ and $\mathcal S_\alpha$ in terms of $\ell_\alpha(X)$. One would similarly expect a lower bounded on the distance between $X$ and $\partial\overline{\T(S)}$ in terms of $\ell_\sys(X)$. Bounds of this type were first obtained by Wu. Before stating Wu's result we define the in-radius of the Teichm\"uller space $\Teich(S)$ by
$$\mbox{\rm InRad}(\Teich(S)) = \max_{X \in \mathcal \Teich(S)} d_{\rm WP}(X,\partial\overline{\T(S)}).$$
Then Wu proves:
\begin{theorem}[Wu, {\cite{Wu:inradius}}]
There exists a universal constant $K$ such that for all $X,Y \in \overline{\T(S)}$ we have
$$\left|\sqrt{\ell_{\rm sys}(X)} - \sqrt{\ell_{\rm sys}(Y)}\right| \le K d_{\rm WP}(X,Y).$$
Therefore
$$d_{\rm WP}(X, \partial\overline{\T(S)}) \ge \frac1K \sqrt{\ell_\sys(X)}$$
and
$$\operatorname{InRad}(S) \ge \frac1K\sqrt{{\rm sys}(S)}.$$
\end{theorem}
By Theorem \ref{wgrad}, for any length function the gradient of $\sqrt{\ell_\alpha}$ is uniformly bounded when the length of the curve is bounded so one would expect a similar statement to hold for $\sqrt{\ell_\sys}$ where the bound depends on $\sys(S)$. What is surprising is that there is a bound independent of topology. 

Here we will show that $\sqrt{\ell_\sys}$ is $1/2$-Lipschitz and we will also give precise asymptotics for $d_{\rm WP}(X, \partial\overline{\T(S)})$ as $\ell_\sys(X) \to \infty$. A key observation in Wu's work is that when a curve is a systole there are improved lower bounds on the width of embedded collars and this leads to better gradient bounds for length functions at systoles. This same observation will be central to our work.

One extra complication is that the systole function is not smooth. However it has enough regularity that we can still discuss its gradient in a modified form that will still satisfy the lower bounds from Lemma \ref{gradflow}. We define
$$\|\nabla \ell_{\rm sys}(X) \| = \underset{\gamma \in {\rm sys}(X)}{\max} \|\nabla\ell_\gamma(X)\|$$
where $\sys(X)$ is the set of curves $\alpha$ that are systoles for $X$. Note that $\sys(X)$ is a finite set so the maximum is always defined.
\begin{lemma}\label{sysflow}
Assume that $U$ is an integrable function with
$$\|\nabla \ell_\sys(X)\| \le U(\ell_\sys(X)).$$
Then for any $X, Y \in \T(S)$ we have
$$d_{\rm WP}(X,Y) \ge \left|\int_{\ell_\sys(X)}^{\ell_\sys(Y)} \frac1{U(s)} ds\right|.$$
\end{lemma}

{\bf Proof:} Let $X_t$ be a smooth path from $X$ to $Y$ parameterized by $[0,1]$ and for each curve $\alpha$ let 
$$f_\alpha(t) = \ell_\alpha(X_t) \qquad \mbox{and} \qquad f_\sys(t) = \ell_\sys(X_t).$$
The path is a compact set in $\T(S)$ and as a smooth function is Lipschitz when restricted to a compact set, each $\ell_\alpha$ will be Lipschitz on the image and therefore each $f_\alpha$ will also be Lipschitz. Furthermore, on a compact set $\ell_\sys$ is the minimum of finitely many length functions so $f_\sys$ is the minimum of finitely many $f_\alpha$. As the minimum of finitely many Lipschitz functions is also Lipschitz we have that $f_\sys$ is Lipschitz. By standard results in analysis $f_\sys$ is differentiable almost everywhere and satisfies the fundamental theorem calculus. Also, as $f_\sys$ is the minimum of finitely many $f_\alpha$ whenever $f'_\sys(t)$ exists we have
$$f'_\sys(t) = f'_\alpha(t) = d\ell_\alpha(\dot X_t)$$
for some $\alpha \in \sys(X_t)$.
Therefore
$$|f'_\sys(t)| \le \|\nabla\ell_\sys(X_t)\| \cdot \|\dot X_t\|$$
when the derivative is defined. The rest of the proof the follows exactly as in Lemma \ref{gradflow}. \eproof

While Theorem \ref{grad} gives bounds on $\|\nabla\ell_{\rm sys}\|$ these bounds can be significantly improved. In particular, for any closed geodesic $\gamma$ on a hyperbolic surface $X$, the collar lemma gives a uniform lower bound on $r_\gamma(X)$ the radius of an embedded collar about $\gamma$ depending only on $\ell_\gamma(X)$. If $\gamma$ is a systole then this radius is bounded below by $\ell_\gamma(X)/4$. For the usual collar lemma the width of the collar decreases to zero as the length grows, in contrast to here where the collar width of the systole limits to infinity. Combining this and Corollary \ref{tgrad} we can improve our upper bounds on the gradient of $\ell_\gamma$ at $X$. We first record the lower bound on the radius of collars of systoles in the following lemma.
\begin{lemma}\label{sys_collar}
Let 
$$r_{\rm sys}(t) = \max\left\{ t/4, \sinh^{-1}\left(\frac{1}{\sinh(t/2)}\right)\right\}.$$
If $\gamma \in {\rm sys}(X)$ then
$$r_\gamma(X) \ge r_{\rm sys}(\ell_\gamma(X)) \ge r_{\rm sys}(L_0)$$
where $L_0$ is the unique positive number with $\sinh(L_0/4)\sinh(L_0/2) = 1$.
\end{lemma}
Combined with Corollary \ref{tgrad} we then have:
\begin{cor}\label{sysgrad}
Let $G$ be the function from Corollary \ref{tgrad}. Then
$$\|\nabla\ell_\sys(X)\|^2 \le \frac{2\ell_\sys(X)}{\pi}(1+G(r_\sys(X)).$$
\end{cor}
Mimicking the definition of the function $H(a,b)$ that we used to bound from below the distance between level sets of lengths functions we define 
$$\qquad H_{\rm sys}(a,b) = \int_a^b \frac{dt}{\sqrt{\frac{2t}{\pi}\left(1+G(r_{\rm sys}(t)\right)}}.$$
As before we further define $H_{\rm sys}(t) = H_{\rm sys}(0,t)$. We also let
$${\mathcal S}^L_\sys = \ell_\sys^{-1}(L) \subset \overline{\T(S)}$$
be the level sets of $\ell_\sys$ and note that ${\mathcal S}^0_\sys = \partial \overline{\T(S)}$.

Note that if $\ell_\sys(X) = b$ and $\gamma \in \sys(X)$ then by Theorem \ref{levelcurve}, for all $a \in [0,\infty)$ we have $d_{\rm WP}(\mathcal S_\gamma^a, X) \le K(a,b)$. As $\ell_\sys \le \ell_\gamma$ if $a\le b$ then, since $\ell_\sys$ is continuous, $d_{\rm WP}({\mathcal S}_\sys^a,X) \le K(a,b)$. In particular, we don't need to modify $K(a,b)$ for the systole function and  we have:
\begin{theorem}\label{syslevelcurve}
If $0\le a< b$ and $X\in \mathcal S_\sys^b$ then
$$H_\sys(a,b) \le d_{\rm WP}({\mathcal S}_\sys^a, {\mathcal S}_\sys^b) \le d_{\rm WP}(\mathcal S_\sys^a,X) \le K(a,b).$$
\end{theorem}
Recall that $K(a,b) = \sqrt{2\pi}\left(\sqrt b - \sqrt a\right)$. It will be useful  to estimate $H_\sys(a,b)$.
\begin{prop}\label{Hestimates}
If $0\le a< b$ then
$$H_\sys(a,b) \ge 2\left(\sqrt b - \sqrt a\right)$$
and
$$\sqrt{\frac2\pi} \le \frac{H_\sys(t)}{\sqrt{2\pi t}}$$
with
$$\lim_{t\to 0} \frac{H_\sys(t)}{\sqrt{2\pi t}} = \lim_{t\to\infty} \frac{H_\sys(t)}{\sqrt{2\pi t}} = 1.$$
\end{prop}

{\bf Proof:} We note that the $r_\sys(t)$ is the maximum of a monotonically increasing and monotonically decreasing function so it is minimized where the two functions agree. That is the minimum of $r_\sys(L_0) = L_0/4$ is the minimum where $L_0$ is the unique positive solution to $\sinh(L_0/4)\sinh(L_0/2) =1$. Therefore
$$H_{\rm sys}(a,b) \ge \sqrt{\frac{2\pi}{1+G(L_0/4)}}\left(\sqrt{b}-\sqrt{a}\right).$$
To evaluate the constant term on right we need to solve $\sinh(L_0/4)\sinh(L_0/2) = 1$ for $L_0$ and the evaluate the function $G$ at $L_0/4$. The function $G$ is an elementary function and can be (rigorously) evaluated using Mathematica to get
$$\sqrt{\frac{2\pi}{1+G(L_0/4)}}\simeq 2.00423$$
and, in particular, it is greater than two. Both inequalities then follow.

For the two limits we observe that $r_\sys(t)$ tends to infinity both as $t\to 0$ and $t\to\infty$ while
$$\lim_{t\to\infty} G(t) = 0.$$
The two limits follow. \eproof

We note that in  \cite[Theorem 1.4]{Wu:inradius} Wu obtains  similar bounds to  Theorem \ref{syslevelcurve}. In both Theorems the upper bound is the same and follows directly from  the lower bound in Riera's formula. In  \cite[Theorem 1.4]{Wu:inradius} the  lower bound is also uniformly comparable to $\sqrt{b}-\sqrt{a}$ as in Theorem \ref{syslevelcurve}.

\begin{figure}[htbp] %  figure placement: here, top, bottom, or page
   \centering
   \includegraphics[width=3in]{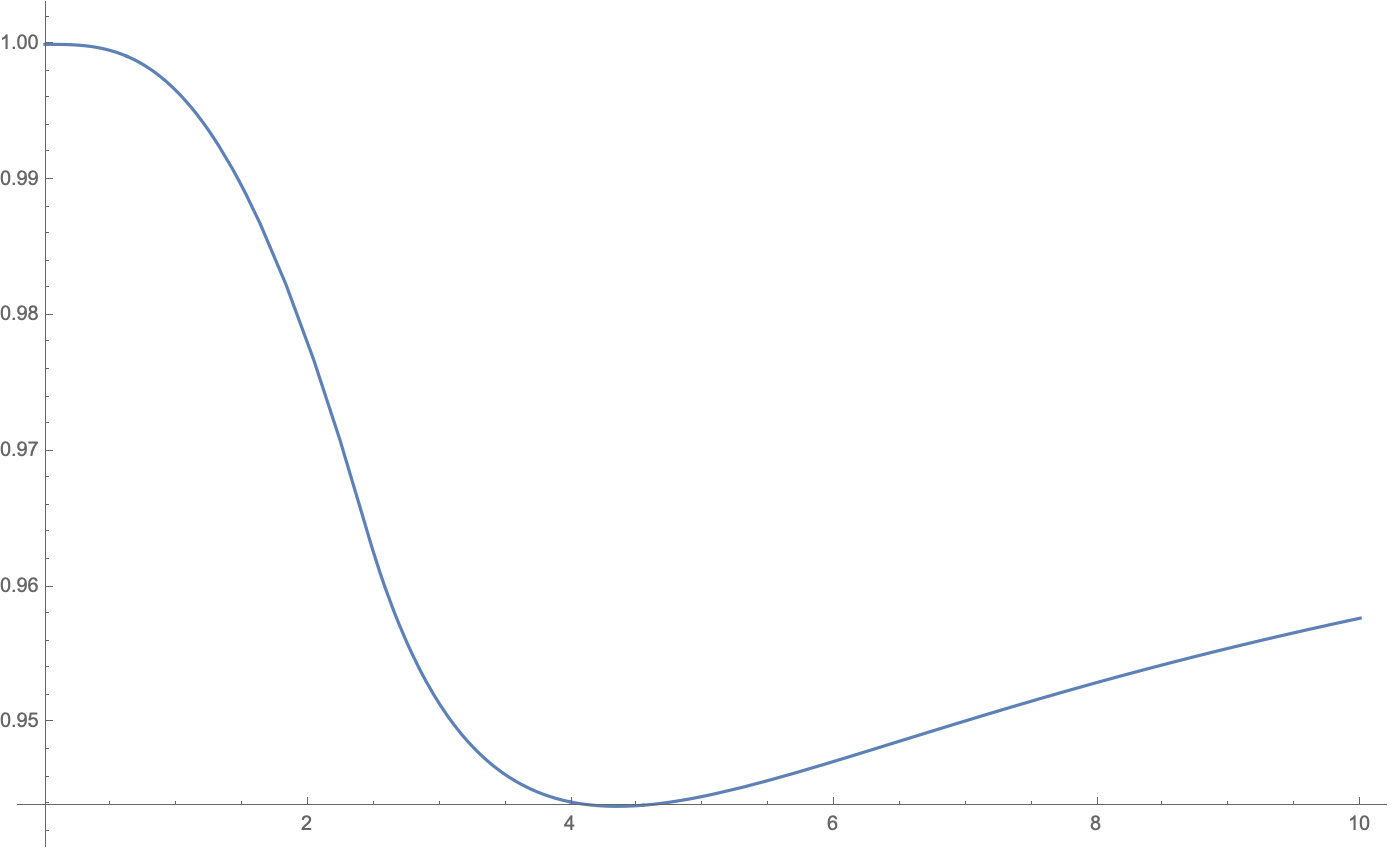} 
   \caption{The graph of $\frac{H_\sys(t)}{\sqrt{2\pi t}}$}
   \label{cfunction}
\end{figure}

{\em Remark:} A more detailed analysis of the function $\frac{H_\sys(t)}{\sqrt{2\pi t}}$  shows that it has a unique critical point which is therefore a global minimum. Evaluating at this minimum gives $\frac{H_\sys(t)}{\sqrt{2\pi t}} \geq .94$ (see Figure \ref{cfunction}).

As an immediate corollary to Theorem \ref{syslevelcurve} we have:
\begin{corollary}\label{lipschitz}
The function $\sqrt{\ell_\sys}$ is 1/2-Lipschitz.
\end{corollary}
We note that using different methods, Wu shows that $\sqrt{\ell_{\rm sys}}$ is Lipschitz with constant $.5492$ for the closed case $S_g$ (see  \cite{Wu:syslip}).

 We also obtain bounds on ${\rm InRad}(S_{g,n})$. For this we apply our work here to bounds on $\sys(S_{g,n})$. For example when $n$ is fixed by \cite{BMP:systole}  we have
  $$\lim_{g\rightarrow \infty} {\rm sys}(S_{g,n}) = \infty.$$  
If $g$ is fixed then $\sys(S_{g,n})$ is uniformly bounded (also see \cite{BMP:systole}). However, it is uniformly bounded below by $2\epsilon_2$. 
 
Thus we have:
\begin{corollary}
For any hyperbolic surface $S$ we have
$$\sqrt{\frac{2}{\pi}}\leq\frac{H_\sys(\sys(S))}{\sqrt{2\pi\sys(S)}} \le \frac{{\rm InRad}(\T(S))}{\sqrt{2\pi\sys(S)}} \le 1,$$
and therefore
$$\lim_{g\rightarrow \infty} \frac{{\rm InRad}(\Teich(S_{g,n}))}{\sqrt{2\pi{\rm sys}(S_{g,n})}} = 1$$
and
$${\rm InRad}(\Teich(S_{g,n})) \geq H_\sys(2\epsilon_2) = H(2\epsilon_2) \simeq 3.27466.$$

\end{corollary}
We note that in \cite[Theorem 1.2]{Wu:inradius} it was shown that ${\rm InRad}(\T(S))$ is uniformly bounded below without producing a concrete bound.
We also remark that that, as in Theorem \ref{smallsep2}, using the fact that we obtain improved lower bounds on the width of collar neighborhoods of separating curves one can show that
$${\rm InRad}(\Teich(S_{0,n})) \geq H_s(4\epsilon_2) \simeq 4.63108.$$

\subsubsection*{Computation}
The calculation of $H, H_s, H_{\rm sys}$ are by numerical integration using Mathematica. The integrand in each  can be written in terms of $F$ where $F(t) = a(T)u(t)v(t)$ where $T = 2\sinh^{-1}(1/\sinh(t/2))$. The functions  $a,u$, and $h$ are elementary functions involving trigonometric, exponential and log functions. To calculate the function $a$ for $t$ small with precision we cannot use its description in terms of basic functions and must instead use a series expansion. The reason for this is that although $a$ is monotonic and $a(0) = 8/3$, the expression for $a$ for small $t$ is the difference of two  large numbers  with the  computation being of the form $(t^{-4} + 8/3\pi^2)-t^{-4}$. To avoid this problem and have arbitrarily high precision, we use the series for the function $\hat a$ introduced in Lemma \ref{eRiera} and the relation 
$$a(T) = \hat a\left(e^{-T}\right) = \hat a\left(\frac{\cosh(t/2)-1}{\cosh(t/2)+1}\right).$$

\begin{figure}[htbp] 
   \centering
   \includegraphics[width=3in]{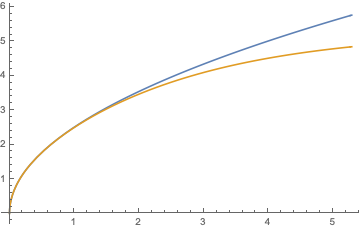} 
   \caption{Graph of H versus K}
   \label{HKgraph}
\end{figure}
See Figure \ref{HKgraph} for a comparison of $H(t)$ and $K(t) = \sqrt{2\pi t}$.

\section*{Appendix: Closed geodesics in the moduli space of the punctured torus}
Our methods can also be used to obtain lower bounds on the minimal Weil-Petersson translation length of a pseudo-Anosov mapping class acting on Teichm\"uller space. We demonstrate the method on the Teichm\"uller space of punctured tori. For a surface of higher complexity the basic idea will still work but it be harder to get explicit estimates.

Let $T$ be the punctured torus and 
$$\psi\colon T\to T$$
a pseudo-Anosov mapping class. By \cite{das:went:wp} there is a unique $\psi$-invariant geodesic $\gamma$ in the Weil-Petersson metric on $\Teich\left(S_{1,1}\right)$. This will descend to a closed geodesic in the moduli space $\mathcal M_{1,1}$. We can use our estimates to give a lower bound on the length of the shortest such geodesic.

We  identify $\Teich(T)$ so that $\psi$ can be represented by an element of $SL_2(\mathbb Z)$:
$$\psi = \left( \begin{array}{cc} a&b\\ c & d\end{array}\right).$$
We can conjugate $\psi$ so that the axis $\gamma$ crosses the imaginary axis at some punctured torus $X$. (This is equivalent to $b/c>0$.) Then $X$ is rectangular: the $(1,0)$-curve and $(0,1)$-curve are represented by geodesics $\alpha$ and $\beta$ that meet orthogonally at a single point. A standard calculation shows that
$$\sinh(\ell_\alpha(X)/2)\sinh(\ell_\beta(X)/2) = 1.$$
One of these two curves will be the systole on $X$ (with the other the second shortest curve).
In fact this is exactly the situation where the collar lemma is optimal: the width of the collar about $\alpha$ is $\ell_\beta(X)$  In a particular if $i(\alpha, \gamma) = k$ then
$$\ell_\gamma(X) \ge k \ell_\beta(X).$$
We have a similar statement when we switch the roles of $\alpha$ and $\beta$.

As $X$ lies on the axis $\gamma$ the translation length of $\psi$ is $d_{\rm WP}(X, \psi(X))$. To bound this distance from below we observe that for any curve $\ell_{\psi(\gamma)}(X) = \ell_\gamma(\psi(X))$. We assume that  $\alpha$  is the shortest curve.

 If  $i(\alpha, \psi(\alpha)) \ge 2$ then 
$$\ell_\alpha(X) \le 2\epsilon_2 \mbox{ and } \ell_{\psi(\alpha)}(X)\ge 2\cdot 2\epsilon_2$$
so by Lemma \ref{levelcurve}
\begin{eqnarray*}
d_{\rm WP}(X, \psi(X)) &\ge& d_{\rm WP}(\mathcal S^{2\epsilon_2}_\alpha, \mathcal S^{4\epsilon_2}_\alpha)\\
& \ge& H(2\epsilon_2,4\epsilon_2)\\
& \ge & 1.06205
\end{eqnarray*}
It follows that  for $\psi$ with $i(\alpha, \psi(\alpha)) \ge 2$ then
$$\|\psi\|_{\rm WP} \geq 1.06205.$$
Otherwise as $\psi(\alpha) = (a,c)$ then $|c| = 1$ and $\psi^2(\alpha) = (a^2+bc, c(a+d))$. As $|a+d| > 2$ then $i(\alpha,\psi^2(\alpha)) = |c(a+d)| =|a+d| \geq 3$ and 
$$\ell_\alpha(X) \le 2\epsilon_2 \mbox{ and } \ell_{\psi(\alpha)}(X)\ge 3\cdot 2\epsilon_2.$$
Therefore 
\begin{eqnarray*}
d_{\rm WP}(X, \psi^2(X)) &\ge& d_{\rm WP}(\mathcal S^{2\epsilon_2}_\alpha, \mathcal S^{6\epsilon_2}_\alpha)\\
& \ge& H(2\epsilon_2,6\epsilon_2)\\
& \ge & 1.56949.
\end{eqnarray*}

Therefore in general
$$\|\psi\|_{\rm WP} \geq \frac{1.56949}{2} \geq  .78474.$$

In \cite{Brock:Bromberg:vol}, the second author and Brock give a  lower bound on the systole for of $\Teich(S_{g,n})$ using renormalized volume and the lower bound for the volume of a hyperbolic 3-manifold. They prove

\begin{theorem}[Brock-Bromberg, {\cite{Brock:Bromberg:vol}}]
Let $\gamma$ be a closed geodesic for the Weil-Petersson metric on moduli space $\mathcal M_{g,n}$ of the surface $S_{g,n}$ with $n >0$. Then
$$\ell_{\rm WP}(\gamma) \geq \frac{4\mathcal V_3}{3\sqrt{\mbox{Area}(S_{g,n})}}$$
where $\mathcal V_3$ is the volume of the regular ideal hyperbolic tetrahedron.\end{theorem}

We note that for $\mathcal M_{1,1}$, the above theorem gives a  bound of $.53724$ and our  bound   is $.78474$. While a more refined analysis could improve this bound, it seems unlikely that these estimates are close to optimal so we do not include them.

\bibliography{bib,math}
\bibliographystyle{math}
\end{document}